\documentstyle[amssymb]{amsart}

\numberwithin{equation}{section}

\newtheorem {Theorem} 			{Theorem}
\newtheorem {varTheorem}                {Theorem}
\newenvironment {Theorem'}
        {\begin{varTheorem}{\hspace{-3.5mm}}{\bf '}{\hspace{3.5mm}}}
        {\end{varTheorem}}
\newtheorem {Lemma}[equation]     	{Lemma}
\newtheorem {Proposition}[equation]	{Proposition}

\newtheorem {Corollary}			{Corollary}

\theoremstyle{definition}
\newtheorem {Definition}[equation]{Definition}
\newtheorem {Remark}[equation]		{Remark}
\newtheorem {Example}[equation]		{Example}

\newcommand     {\comment}[1]   {}
\newcommand     {\mute}[2] {}
\newcommand     {\printname}[1] {}

\newcommand{\labell}[1] {\label{#1}\printname{#1}}

\def	\ch	{\text{ch}}
\def	\red	{{\operatorname{red}}}
\def	\cut 	{{\operatorname{cut}}}
\def    \coker	{\operatorname{coker}}

\def	\Euler	{\operatorname{Euler}}
\def	\Trace	{\operatorname{Trace}}

\def	\inv	{^{-1}}
\def	\to	{\longrightarrow}
\def	\half	{{\frac12}}
\def	\ssminus	{\smallsetminus}

\def	\C	{{\Bbb C}}
\def	\R	{{\Bbb R}}
\def	\Z	{{\Bbb Z}}
\def	\CP	{{\Bbb C}{\Bbb P}}

\def	\cO	{{\cal O}}
\def	\cJ	{{\cal J}}
\def	\cS	{{\cal S}}

\def	\fg	{{\frak g}}

\def	\spinc	{{\operatorname{spin}^c}}
\def	\spin	{\operatorname{spin}}
\def	\SO	{\operatorname{SO}}
\def	\cF	{{\cal F}}

\def	\Lie	{\operatorname{Lie}}

\def	\AS	{\operatorname{AS}}
\def	\Aroof	{{\hat A}}
\def	\Re	{\operatorname{Re}}
\def	\Im	{\operatorname{Im}}

\def	\tc	{\tilde{c}}
\def	\ti	{\tilde{i}}
\def	\tiplus	{\tilde{i}_+}

\def	\tPhi	{\tilde{\Phi}}
\def	\tZ	{\tilde{Z}}

\def	\bs	{\boldsymbol}
\def	\pr	{\operatorname{pr}}

\newcommand {\dd}[1]	{{\partial \over \partial {#1}}}

\begin{document}


\title[Presymplectic quantization]
{Quantization of Presymplectic Manifolds and Circle Actions}

\author{Ana Canas da Silva}
\author{Yael Karshon}
\author{Susan Tolman}
\thanks{A. Canas da Silva was partially supported by a NATO fellowship.
Her research at MSRI was supported in part by NSF grant DMS 9022140.
Y. Karshon was partially supported by NSF grant DMS 9404404.
S. Tolman was partially supported by an NSF postdoctoral fellowship.}
\thanks{dg-ga/9705008}

\address{Department of Mathematics, University of California,
Berkeley, CA 94720}
\email{acannas@@math.berkeley.edu}

\address{Inst. of Math.,
The Hebrew University, Jerusalem 91904, Israel}
\email{karshon@@math.huji.ac.il}

\address{Department of Mathematics, Princeton University,
Princton, NJ 08544-1000}
\email{tolman@@math.princeton.edu}

\thanks{May 4, 1997}

\begin{abstract}
We prove several versions of ``quantization commutes with reduction''
for circle actions on manifolds that are not symplectic.
Instead, these manifolds possess a weaker structure,
such as a spin$^c$ structure. Our theorems work whenever
the quantization data and the reduction data are compatible;
this condition always holds if we start from a presymplectic
(in particular, symplectic) manifold.
\end{abstract}

\maketitle


\section{Introduction}
\labell{sec:intro}

Morally, {\bf quantization} is a functor that
associates a Hilbert space to every symplectic manifold and linear operators
to functions, subject to certain conditions suggested by physics.
In particular, quantization should associate to every Hamiltonian
action of a group  on a symplectic manifold
a unitary representation of that group.
Symplectic {\bf reduction\/} models the fact that symmetries of a
Hamiltonian dynamical system enable one to reduce the number
of degrees of freedom of the system; it produces a
lower dimensional symplectic manifold.
In the last few years there has been much progress
around a ``conjecture'' of Guillemin and Sternberg~\cite{GS} that
can be loosely phrased as {\bf quantization commutes with reduction}.
This means that the quantization of the reduced space ought to be naturally
isomorphic to the space of fixed vectors in the quantization
of the original manifold.
In this paper we show that the assertion ``quantization commutes
with reduction'' makes sense and is true in a much more general context
than that of symplectic manifolds.
To make this statement precise, both quantization and reduction must
be rigorously defined.
To make this statement true, some compatibility relations
between the data for quantization and the data for reduction
must be assumed.

Although it is impossible to satisfy all the requirements
of quantization,  some recipes are partially successful.
For example, let a group $G$ act in a Hamiltonian fashion
on a K\"ahler manifold $(M,\omega)$.
Suppose that $\omega$ is integral, and let $L$ be a $G$-equivariant
holomorphic line bundle with first Chern class $[\omega]$.
The space of global holomorphic sections of $L$ is a representation of
$G$; this is called the {\bf K\"ahler quantization}.

The first natural extension of this idea, suggested by Bott~\cite{B},
is to replace the space of holomorphic sections by the alternating sum
of sheaf cohomology groups,
$
   \sum_i (-1)^i H^i(M,{\cal O}_L)\ ,
$
where ${\cal O}_L$ is the sheaf of germs of holomorphic sections of $L$.
By Kodaira's theorem,
$H^i(M,{\cal O}_{L^n}) = 0$  for all $i \geq 1$ for sufficiently
large $n$, and we retrieve K\"ahler quantization in the semi-classical
limit.
This alternating sum is a {\bf virtual representation} of $G$,
{\em i.e.\/}, a formal difference of two linear representations.
By Dolbeault's theorem, it is the same as the index
of the Dolbeault $\overline{\partial}$ operator with coefficients in $L$.

A second generalization, that we will call
{\bf almost-K\"ahler quantization}\footnote
{Many authors call this $\spinc$ quantization, but we prefer to reserve
that name for the quantization defined on all $\spinc$ structure.},
applies to symplectic manifolds
$(M,\omega)$ on which a group $G$ acts in a Hamiltonian fashion.
This uses the existence, and uniqueness up to homotopy,  of an
almost complex structure $J$ that is {\bf compatible} with $\omega$,
i.e., such that the bilinear form
$$
   g_p (v,w) = \omega_p (J_pv,w)\ ; \quad v,w \in T_pM
$$
is symmetric and positive definite for all $p \in M$.
(When $J$ is actually a complex structure, compatibility is
equivalant to $\omega$ being a K\"ahler form.) If $\omega$ is integral,
there is a complex line bundle $L$ with first Chern class $[\omega]$.
Choosing additional auxiliary data, this determines an elliptic
operator, ${\cal D}$, whose index only depends on $\omega$ and $J$.
(Bott calls it the ``rolled-up Dolbeault operator'', whereas
Duistermaat coined the term ``Dolbeault-Dirac operator''.)
Almost-K\"ahler quantization is defined to be the index of this
operator:
$$
 [\ker \cal D] - [\coker \cal D]\ .
$$
The kernel and cokernel of $\cal D$ are linear
$G$-representations.\footnote{Recently, Borthwick and Uribe~\cite{BU}
have proved vanishing theorems for the cokernels, showing that
almost-K\"ahler quantization yields honest vector spaces semiclassically.}
The index, as a virtual representation,
is determined by its character function, $\chi : G \to \C$, given by
$\chi(g) = \Trace(g|_{\ker{\cal D}}) - \Trace(g|_{\coker{\cal D}})$.
The multiplicity of $a$ in the index
is its multiplicity in $\ker \cal D$ minus its multiplicity in $\coker \cal D$.
Guillemin and Sternberg observed that the Riemann-Roch number,
which computes this index, is in fact
a symplectic invariant, and hence can be used to quantize a symplectic
manifold. See Vergne~\cite{V1} and Guillemin~\cite{G}.

One can also define quantization as the index of other differential operators,
such as the Dirac spin or spin$^c$ operators.
Bott used the Dirac spin$^c$ operator in this context
back in the 60's \cite{Bott-spinc}.
Representation theorists, attempting to construct irreducible representations
by quantizing coadjoint orbits, used the Dirac spin operator in the 1970's;
see, {\em e.g.\/}, \cite{par} and references therein.
Duflo prefered  the spin$^c$ operator for this purpose.

We begin our paper by discussing, in section \ref{sec:quantization},
different possible extensions  of symplectic
quantization to non-symplectic settings.
We work with each of the following
{\bf quantization data\/} (in increasing generality):
\begin{enumerate}
\item
an almost complex structure and a complex line bundle;
\item
a stable complex structure, an orientation, and a complex line bundle;
or
\item
a spin$^c$ structure.
\end{enumerate}
Each such data determines a Dirac operator, $\cal D$, which acts
on the sections of certain vector bundles over $M$.
The corresponding quantization is defined to be the index of
this operator.
The definition of the Dirac operator requires some additional structure
(connection, metric), but its index is independent of these choices.

Let a compact Lie  group $G$ act on a symplectic manifold $(M,\omega)$
with a moment map $\Phi : M \to \fg^*$, where $\fg^*$ is the
dual to the Lie algebra of $G$.
For every regular value $a \in \fg^*$, the {\bf reduced space},
$M^a_\red := \Phi^{-1}(G \cdot a)/G$, is naturally a symplectic
orbifold, and it is a symplectic manifold when $G$ acts freely
\cite{MW}.
In section \ref{sec:reduction}
we show how $S^1$-equivariant quantization data
descends to quantization data on the reduced space.
In the absence of a moment map and its level sets,
we define the reduced space to be
the quotient by the circle action of an arbitrary invariant hypersurface
that splits our manifold into two pieces.

Guillemin and Sternberg proved that K\"ahler quantization commutes
with reduction 
in the same paper where they stated their conjecture~\cite{GS}.
More specifically, they proved that
the dimension of the quantization of the reduced space at $a=0$
is equal to the  multiplicity of $0$ in the quantization
of the original manifold.
Recently, much progress has been made in proving that
almost-K\"ahler quantization commutes with symplectic reduction
under various assumptions;
Duistermaat, Guillemin, Jeffrey, Kirwan, Meinrenken, Sjamaar,
Vergne and Wu
(\cite{M1}, \cite{4}, \cite{G}, \cite{JK}, \cite{V2}, \cite{M2},
\cite{MS})
covered the cases where one starts with a symplectic orbifold,
where the level for reduction is singular, and where
the group acting is an arbitrary compact Lie group.
A survey of these developments, as well as an excellent introduction
to the problem, can be found in \cite{S}.

We extend this theorem  by proving that our three versions
of quantization commute with our generalized notion of
reduction as long as
certain components of the fixed point set  lie
on the correct side of the splitting hypersurface.
The precise condition depends on the quantization data;
see section \ref{sec:main}.
Our main theorem, that $\spinc$ quantization commutes
with reduction, implies the corresponding theorems for stable complex
and almost complex quantization.
Hence it implies the standard theorem for almost-K\"ahler
quantization.\footnote{
Note that other authors use the term ``$\spinc$ quantization''
to refer only to this very special case.}

Each of these three types of
quantization data gives a {\bf presymplectic\/} form
(that is, a closed two-form)\footnote{
We do not insist, as some authors do, that the closed two-form
have constant rank.}
on our manifold; more precisely,  each gives the cohomology class
represented by such a form.  We think
of the quantization as a quantization of the presymplectic manifold
and of the quantization data as a ``polarization'' -- auxiliary
data used to define the quantization.
Although the presymplectic forms are not essential for quantization
nor reduction nor ``quantization commutes with reduction'',
they are relevant
because traditionally quantization and reduction only applied to
symplectic manifolds, and because presymplectic forms
provide interesting examples.

Presymplectic forms come up naturally in many mathematical contexts:
as curvature forms of line bundles,
as pull-backs of K\"ahler forms via complex blow-ups,
on connected sums of symplectic manifolds,
as the averages of symplectic forms by group actions \cite{GK1},
as folding type structures \cite{CGW},
and in the framework of symplectic cobordism \cite{GGK}.
For an application to representation theory see \cite{GK1}.
Presymplectic forms also come up in physics \cite{constraints}.

To a group action on a presymplectic manifold
one can associate a moment map, just as one does on a symplectic manifold.
The zero level set of the moment map is a
hypersurface that automatically satisfies the compatibility conditions
mentioned earlier. Thus ``presymplectic quantization commutes with
reduction".
We provide the details in section~\ref{sec:corollaries}.
This theorem was proven in the special case of almost-K\"ahler
quantization of toric manifolds in~\cite{KT}, then generalized
to completely integrable torus actions and $\spinc$ quantization
in~\cite{GK2}.

In sections~\ref{sec:cutting} and~\ref{proof-spinc} we prove
our main theorem, that $\spinc$ quantization commutes with reduction.
Our proof is modeled on
the simple and elegant proof of Duistermaat-Guillemin-Meinrenken-Wu \cite{4}
for the symplectic case.
The result follows from applying the index theorem of Atiyah-Segal-Singer
to the original space, to the ``Lerman cut space''
(see section~\ref{sec:cutting}), and to the reduced
space, and comparing terms.

Spin$^c$ quantization is our favorite, for several reasons:
spin$^c$ quantization of a presymplectic manifold
depends on no auxiliary data; it only depends on the
manifold, the orientation, and the cohomology class of the two-form.
The theorem that ``quantization commutes with reduction"
is the strongest in the $\spinc$ case; as mentioned earlier,
it implies similar theorems for almost complex and stable complex
quantization. For presymplectic manifolds, the $\spinc$ theorem
is the most robust in that it gives the biggest freedom to choose
the level set at which we can reduce; see Remarks~\ref{vary0} and~\ref{vary1}.

We have designed this paper so that the reader who is only interested in
the $\spinc$ case does not have to read about the other cases.  To read the
paper this way, simply skip all the subsections with the word ``complex''
in the title, also skip Theorems \ref{theorem-a-c} and \ref{theorem-s-c},
the remarks immediately after them,
and their proofs in section \ref{sec:main},
and additionally skip Corollaries \ref{coroll1} and \ref{coroll2}
in section \ref{sec:corollaries}.

In this paper we restrict to circle actions on manifolds and insist that
the circle act freely on the reducible hypersurface, so that the reduced
space is a smooth manifold.  However, we expect
our results to generalize to torus actions and to the case
where the  original space is an orbifold or the action is only
locally free.

As for actions of nonabelian groups,
the most straightforward generalization of the theorem
is true in the symplectic case \cite{M2} but is false in the
almost complex, stable complex, and spin$^c$ presymplectic cases;
counterexamples have been found by Vergne and by Vogan,
for the action of $SU(2)$ on $\CP^1$.
(See \cite{MS}.)
The general presymplectic nonabelian case was recently resolved
by Tolman \cite{tolman:nonabelian}.

\subsection*{Acknowledgements}

The authors wish to acknowledge their gratitude for the
useful discussions, suggestions, and support of Victor Guillemin,
Viktor L. Ginzburg, and Eckhard Meinrenken.


\section{Quantization}
\labell{sec:quantization}

In this section we describe several ways to extend almost-K\"ahler
quantization
(see section~\ref{sec:intro}) and to quantize non-symplectic manifolds
via indices of elliptic differential operators.
This is mostly review of known material, but is necessary for the
rest of the paper. In \S\ref{presymplectic quantization} we interpret
these as quantizations of presymplectic manifolds.


\subsection{Almost complex quantization}
\labell{almost complex quantization}

Let $M$ be a compact manifold, $L \to M$ a complex line
bundle, and $J : TM \to TM$ an almost complex structure.
Choose a metric on $M$, a Hermitian structure on $L$,
and a Hermitian connection on $L$.
This data determines an elliptic operator, $\cal D$,
called the {\bf twisted Dolbeault-Dirac operator}
(see \cite[\S 9]{GK2} or \cite{D}).
The {\bf almost complex quantization} of $(M,L,J)$
is the index of this operator:
$$
	Q(M,L,J) := [\ker{\cal D}] - [\coker{\cal D}] \ .
$$
If the above data is equivariant with respect to an action of
a Lie group $G$, the index is a  virtual representation of $G$.
The index is independent of the auxiliary choices of metric
and connection.


\subsection{Spin$^c$ quantization}
\labell{spinc quantization}

The group $\spin(n)$ is the connected double cover of the group $\SO(n)$.
Let $q:\spin(n) \to \SO(n)$ be the covering map, and let $\{1,\epsilon\}$
be its kernel. The group $\spinc(n)$ is the quotient
$(\spin(n) \times U(1))/\Z_2$,
where $\Z_2$ is generated by $(\epsilon,-1)$.
There are two canonical group homomorphisms:
$$
   \begin{array}{crclccrcl}
         f: & \spinc(n) & \to & \SO(n) & \;\;\; \mbox{and}  \;\;\; &
         \det: & \spinc(n) & \to & U(1) \\
         & [s,u] & \mapsto & q(s) & &
         & [s,u] & \mapsto & u^2
   \end{array}
$$
(the square brackets indicate class modulo
 $\Z_2$).

Let $V \to M$ be a real oriented vector bundle of rank $n$,
equipped with a fiber metric,\footnote{
It would be more elegant, but less standard,
to follow Duflo and Vergne in using {\bf quantum line bundles}
(see \cite{V2}) instead of spin$^c$ structures.
These eliminate the unnecessary choice of a metric.
}
and let $\cF(V)$ be its oriented orthonormal frame bundle.
The group $\spinc(n)$ acts on $\cF(V)$ from the right
via the homomorphism $f$.
A {\bf $\bs \spinc$ structure} on the vector bundle $V$ is a
$\spinc(n)$-principal bundle, $P \to M$, together with a
$\spinc(n)$-equivariant bundle map, $p : P \to \cF(V)$.
The line bundle $L = P \times_{\det} \C$
is called the {\bf line bundle associated to $\bs P$}.
A spin$^c$ structure on an oriented Riemannian manifold
is a $\spinc$ structure on its tangent bundle.

If a group $G$ acts on the bundle $V$, preserving the fiber metric,
we define an {\bf equivariant $\bs \spinc$ structure} on $V$ by
requiring $P$ to be a $G$-equivariant principal bundle and the
map $p$ to be $(G \times \spinc(n))$-equivariant.

A $\spinc$ structure $(P,p)$ on a compact oriented Riemannian
manifold $M$, plus a
connection on $P$, determines an elliptic operator, $\cal D$,
called the {\bf  $\bs \spinc$-Dirac operator}
(see \cite[\S 9]{GK2} or \cite{D}).
The {\bf $\bs \spinc$ quantization} of $(M,P,p)$ is
the index of this operator:
\begin{equation} \labell{spinc index}
	Q(M,P,p) :=  [\ker {\cal D}] - [\coker {\cal D}]
\end{equation}
If $(P,p)$ is a $G$-equivariant spin$^c$ structure,
and the connection on $P$ is $G$-invariant,
the index is a virtual $G$-representation.
The index is independent of the choice of connection on $P$.

\begin{Lemma} \labell{a-c implies spinc}
Let $M$ be an oriented Riemannian manifold with a $G$ action.
There is an equivariant $\spinc$ structure on $M$ naturally associated
to  every invariant almost complex structure on $M$ and equivariant
line bundle $L$ over $M$.
The line bundle associated to this $\spinc$ structure
is $L \otimes L \otimes K^*$, where $K$ is the canonical
line bundle, $\bigwedge^n_\C T^*M$.
The $\spinc$ quantization then coincides with the
almost complex quantization.
\end{Lemma}

\begin{pf}
Apply the standard construction
(see~\cite{LM}, appendix D, example D.6)
to the tangent bundle $TM$.
The quantizations are the same because the Dolbeault-Dirac
operator and the $\spinc$-Dirac operator $\cal D$ have the same
principal symbol and hence the same index; see~\cite[\S 5]{D}.
\end{pf}

\begin{Remark} \labell{global sign}
There are two different elliptic operators
on an even dimensional $\spinc$ manifold: $D^+:S^+ \to S^-$ and
$D^- = (D^+)^* : S^- \to S^+$, where $S^+$ and $S^-$ are the so-called positive
and negative spinors.
In dimensions 2 mod 4, there is no universally recognized standard as to
which half of the spinor bundle should be called $S^+$.
We follow the conventions of Lawson-Michelson~\cite{LM},
Duistermaat~\cite{D}, and the majority of recent authors.
Unfortunately, these disagree with Atiyah-Singer~\cite{ASIII}
and Vergne~\cite{V2}, and thus our formulas
will differ from those in~\cite{ASIII} by $(-1)^n$ in dimension $2n$.
\end{Remark}

We will use the following lemma twice:
to define the quantization of a stable complex structure
in \S\ref{stable complex quantization}, and to define the
reduction of a $\spinc$ structure in \S\ref{sec:reduction}.

\begin{Lemma}
\labell{reduced-spin-c}
Let $V \to M$ be a vector bundle
equipped with a fiber metric and orientation
and with a left action of a Lie group $G$.
Every equivariant spin$^c$ structure on
a direct sum of $V$ with an equivariantly trivial bundle naturally
induces an equivariant spin$^c$ structure on $V$, with the same
associated line bundle.
\end{Lemma}

\begin{pf}
Let $n = \text{rank} (V).$
The splitting $\R^{n+k} = \R^n \times \R^k$ induces
the following commuting diagrams of group homomorphisms:
\begin{equation}
\labell{fjcommutes}
   \begin{array}{ccc}
   \spinc(n) & \hookrightarrow & \spinc(n+k) \\
   {\scriptstyle f} \downarrow &  & {\scriptstyle f} \downarrow \\
   \SO(n) & \hookrightarrow & \SO(n+k)
   \end{array}
\end{equation}
and
\begin{equation}
\labell{detjcommutes}
   \begin{array}{ccc}
   \spinc(n) & \hookrightarrow & \spinc(n+k) \\
   {\scriptstyle \det} \downarrow \phantom{\scriptstyle \det}  &  &
   {\scriptstyle \det} \downarrow \phantom{\scriptstyle \det}      \\
   U(1) & \stackrel{=}{\rightarrow} & U(1)
   \end{array}
\end{equation}

Let $(P',p')$ be an equivariant spin$^c$ structure on $V \oplus \R^k$.
The bundle $P$, which is defined by the following pull-back diagram,
$$
   \begin{array}{ccc}
   P & \to & P' \\
   \phantom{\scriptstyle p} \downarrow {\scriptstyle p}
   & &
   \phantom{\scriptstyle p'} \downarrow {\scriptstyle p'}  \\
   \cF(V) & \stackrel{j}{\hookrightarrow} & \cF(V \oplus \R^k)
   \end{array}
$$
is a $G$-equivariant spin$^c$(n)-principal bundle over $M$.
The map $p := j^* p'$ is $(G \times \spin^c(n))$-equivariant
by \eqref{fjcommutes}.
By \eqref{detjcommutes}, the line bundle associated to $P$
is the same as that associated to $P'$.
\end{pf}


\subsection{Stable complex quantization}
\labell{stable complex quantization}

Stable complex structures are a generalization of
almost complex structures which behave much better under reduction;
see \S~\ref{sec:reduction}.

A {\bf stable complex structure} on a manifold $M$ is an
equivalence class of complex structures on the vector bundles
$TM \oplus \R^k$; two such complex structures, on bundles
$E_1 = TM \oplus \R^{k_1}$ and $E_2 = TM \oplus \R^{k_2}$,
are equivalent if there exist $r_1$, $r_2$ such that
$E_1 \oplus \C^{r_1}$ and $E_2 \oplus \C^{r_2}$
are isomorphic  complex vector bundles.
When a group acts on $M$ we define an
{\bf equivariant stable complex structure\/} by requiring the complex
structures to be invariant and the isomorphism to be equivariant.
Here the group acts on $TM$ by the natural lifting of its action on
$M$, and it acts trivially on the trivial bundles $\R^k$ and $\C^r$.
The {\bf canonical line bundle} of a stable complex structure
is the top wedge $K = \bigwedge^{m}_\C E^*$,
where $E = TM \oplus \R^k$ is a complex vector bundle
of rank $m$ representing the stable complex structure.

An almost complex structure determines
a stable complex structure and an orientation.
A stable complex structure does not determine an orientation.\footnote{
For instance, the almost complex structures $J$ and $-J$
on $\C$ define the same stable complex structure but opposite
orientations.}

\begin{Lemma} \labell{s-c implies spinc}
Let $M$ be an oriented Riemannian manifold with a $G$ action.
There is an equivariant $\spinc$ structure on $M$ naturally associated
to  every invariant stable complex structure on $M$ and equivariant
line bundle $L$ over $M$.
The line bundle associated to this $\spinc$ structure
is $L \otimes L \otimes K^*$, where $K$ is the canonical
line bundle of the stable complex structure.
\end{Lemma}

\begin{pf}
The standard construction
(see \cite{LM}, appendix D, example D.6),
applied to a complex vector bundle that represents the stable complex
structure, yields a $\spinc$ structure on the vector bundle $TM \oplus
\R^k$ with associated line bundle $L \otimes L  \otimes K^*$.
By Lemma~\ref{reduced-spin-c}, this spin$^c$ structure
determines a $\spinc$ structure on $TM$ with the same associated line
bundle.  Different representatives of the stable complex class yield
isomorphic spin$^c$ structures by naturality of the standard
construction in~\cite{LM} and of the pull-back diagram in the proof
of Lemma~\ref{reduced-spin-c}.
\end{pf}

\begin{Remark}
In passing from an almost complex or stable complex structure and a
line bundle to the corresponding spin$^c$ structure we lose
information.
For example,
on $\CP^1 \times \CP^1$
the line bundle $O(m) \boxtimes O(m)$ with the
standard complex structure, and the line bundle
$O(m+2) \boxtimes O(m+2)$ with the opposite complex structure
induce the same $\spinc$ structure.
\end{Remark}

Let $M$ be a compact manifold with an action of a compact Lie group $G$,
$L$ an equivariant complex line bundle over $M$,
$\cJ$ an equivariant stable complex structure on $M$,
and $\cO$ an orientation on $M$.
Let $(P,p)$ be the corresponding spin$^c$ structure
described in Lemma~\ref{s-c implies spinc}
(for an arbitrary Riemannian metric on $M$),
and let $\cal D$ be the corresponding spin$^c$-Dirac operator
(for an arbitrary connection on $P$).
{\bf Stable complex quantization} is the index of this operator,
{\em i.e.\/}, is the virtual $G$-representation:
$$
	Q(M,L,\cJ,\cO) = [\ker {\cal D}] - [\coker {\cal D}].
$$

\begin{Remark}
\labell{a-c-q = spinc-q}
Let $J$ be an {\em almost\/} complex structure on $M$ and $L$ a
complex line bundle over $M$.  Let $\cJ$ and $\cO$ be the
stable complex structure and the orientation induced by $J$.
The $\spinc$ structure associated to the data $(M,L,J)$
by Lemma \ref{a-c implies spinc} coincides with the
$\spinc$ structure associated to the data $(M,L,\cJ,\cO)$
by Lemma \ref{s-c implies spinc}.
Consequently, the almost complex quantization of $(M,L,J)$
is the same as the stable complex quantization of $(M,L,\cJ,\cO)$.
\end{Remark}


\subsection{Quantization of presymplectic manifolds}
\labell{presymplectic quantization}

A {\bf presymplectic manifold} is a manifold $M$
equipped with a {\bf presymplectic form} $\omega$,
{\em i.e.\/}, a closed, possibly degenerate, two-form.
Let the circle group, $S^1$, act on $M$ and preserve $\omega$.
A {\bf moment map} is a map
$\Phi : M \to \R$ such that $d\Phi =  \iota (\eta_M) \omega$,
where $\eta$ the positive primitive lattice element in $\Lie(S^1)$,
and $\eta_M$ is the vector field on $M$ which generates the action
of ${\rm exp} (t\eta)$, for $t \in \R$
(See Appendix~\ref{conventions} for an explanation of our conventions.)

Given a moment map $\Phi$, we can define the quantization of $(M,\omega,\Phi)$
using each of the three types of quantization above.
For this we need to impose an integrality-type condition
on $\omega$ and $\Phi$; the precise condition
depends on whether we perform (almost or stable) complex quantization
or $\spinc$ quantization.


\subsubsection*{Almost complex and stable complex quantizations}

The triple $(M,\omega,\Phi)$ is {\bf complex prequantizable}
if the equivariant cohomology class\footnote{
Here, $u$ is a formal variable.  For definitions of equivariant forms
and equivariant cohomology, please see \cite{BGV}.
}
$[\omega + u\Phi]$ is integral.
This condition precisely guarantees that we can find an equivariant
line bundle $L$ whose equivariant first Chern class is
$[\omega + u\Phi]$.
(The non-equivariant case is explained well in \cite{kostant}.
By applying that to $M \times_{S^1} ES^1$, we get the equivariant case.)

Given an invariant almost complex structure $J$ on $M$, the
almost complex quantization of the presymplectic manifold is
\begin{equation}
\labell{def1}
	Q(M,\omega,\Phi,J) := Q(M,L,J)\ .
\end{equation}
Given an equivariant stable complex structure $\cJ$ on $M$
and an orientation $\cO$ on $M$,
the stable complex quantization of the presymplectic manifold
is
\begin{equation}
\labell{def2}
	Q(M,\omega,\Phi,\cJ,\cO) := Q(M,L,\cJ,\cO)\ .
\end{equation}
Although there may exist different line bundles
with the same equivariant first Chern class,
the quantizations \eqref{def1} and \eqref{def2}
are independent of this choice.
This follows from the Atiyah-Segal-Singer formula for the equivariant index;
see \S\ref{proof-spinc}.


\subsubsection*{Spin$^c$ quantization}

The triple $(M,\omega,\Phi)$ is {\bf spin$^{\mathbf c}$ prequantizable}
if the equivariant cohomology class $2 [\omega + u \Phi]$ is integral
and the modulo 2 reduction of $2 [\omega + u \Phi]$ is
the second equivariant Stiefel-Whitney class of $M$.
For any orientation $\cO$ and Riemannian metric on $M$,
this condition precisely guarantees that
there exists a $G$-equivariant spin$^c$ structure $(P,p)$ on $M$
such that the equivariant  first Chern class of its associated line bundle is
$2 [\omega + u \Phi]$; see \cite{LM}.
The $\spinc$ quantization of the presymplectic manifold is defined to be
\begin{equation}
\labell{def3}
        Q(M,\omega,\Phi,\cO) := Q(M,P,p).
\end{equation}
It is independent of the
the metric on $M$ and  the $\spinc$ structure $(P,p)$.
This follows from the Atiyah-Segal-Singer formula
for the equivariant index; see \S\ref{proof-spinc}.

\begin{Remark}
\labell{shift}
Let $S^1$ act on $M$, and  consider a $\spinc$ prequantizable
equivariant two-form,  $\omega + u \Phi$.
For any weight $a \in \Z$, the equivariant  two-form
$\omega + u (\Phi- a)$ is also  $\spinc$ prequantizable;
its quantization differs from that of $(M,\omega,\Phi,\cO)$
by tensoring with the trivial $S^1$-representation with weight $a$. Hence
the multiplicity of $a$ in $Q(M, \omega, \Phi, \cO)$ is equal to the
multiplicity of $0$ in $Q(M,\omega,\Phi-a, \cO)$.
This fact is needed to make the shifting trick work, and hence to make
quantization commute with reduction at all weights; see
section~\ref{sec:corollaries}.
Note that this remark would be false without the factor of $2$
in the definition of spin$^c$ quantization.

An analogous remark holds for almost and stable complex structures,
that is, the multiplicity of $a$ in $Q(M,\omega,\Phi,J)$ is
the multiplicity of $0$ in $Q(M,\omega,\Phi-a,J)$, etc.
\end{Remark}

\begin{Remark}
A line bundle $L$ and  an almost (or stable) complex structure $J$ on $M$
naturally determine a $\spinc$
structure $(P,p)$ such that the quantizations
$Q(M,L,J)$ and $Q(M,P,p)$  coincide.
Nevertheless, we interpret these as the
quantizations of two {\em different\/} presymplectic forms.
Specifically, $Q(M,L,J)$ is the quantization of a two-form representing
the first Chern class of $L$,
whereas $Q(M,P,p)$ is the quantization of a two-form representing the
first Chern class of
$L$ {\em plus} half the first Chern class of the dual to the canonical bundle
of the almost complex structure.
(See Lemma~\ref{a-c implies spinc} and Example~\ref{CP2}.)

\end{Remark}

\begin{Example} \labell{CP2}
Let $[\sigma] \in H^2(\CP^2)$ generate the integral cohomology.
The dual to the canonical bundle has first Chern class $3[\sigma]$.
A presymplectic form
$\omega$ is spin$^c$-quantizable if and only if
$[\omega] = (k + 3/2) [\sigma]$ for some $k \in \Z$;
in particular, the $\spinc$-quantizable two-forms are {\em not\/} integral.
The spin$^c$ quantization of  $(k + 3/2) \sigma$ coincides with
the almost complex quantization of  $k\sigma$.
\end{Example}


\section{Reduction}
\labell{sec:reduction}

Let the circle act on a manifold $M$.
When $M$ possesses a symplectic form that is invariant under
the action and has a moment map $\Phi: M \to \R$, symplectic reduction
produces a symplectic structure
on the reduced space $M^a_\red = \Phi\inv(a) / {S^1}$ for every
regular value $a$ (see \cite{MW}).
We generalize this idea in two ways:
first, we take the quotient of arbitrary co-oriented hypersurfaces
on which the circle acts freely.
Second, we describe how additional equivariant structures on $M$
descend to this  subquotient of $M$.

\begin{Definition} \labell{def-reducible}
A {\bf reducible hypersurface} in $M$ is a co-oriented submanifold
$Z$ of codimension one that is invariant under the $S^1$ action
and on which this action is free.
The {\bf reduction of $\bs M$ at $\bs Z$}
is the quotient of $Z$ by the circle action, $M_\red := Z / S^1$.
\end{Definition}

It is useful to keep in mind the diagram
\begin{equation} \labell{reduction}
\begin{array}{ccc}
 	Z & \stackrel{i}{\hookrightarrow} & M \\
	{\scriptstyle \pi} \downarrow \phantom{\scriptstyle \pi}&&\\
	M_\red	& &
\end{array}
\end{equation}
where $i$ and $\pi$ are the inclusion and quotient maps, respectively.

\begin{Definition}
\labell{def-splits}
A reducible hypersurface $Z$ is {\bf splitting} if its
complement, $M \ssminus Z$, is a disjoint union of two (not necessarily
connected) open pieces,
$M_+$ and $M_-$, such that positive normal vectors to $Z$ point
into $M_+$ and negative normal vectors point into $M_-$.
We then say that $Z$ {\bf splits} $M$ into $M_+$ and $M_-$.
\end{Definition}

\begin{Example}
\labell{ex-phi}
Let $\Phi: M \to \R$ be a smooth $S^1$-invariant function.
Assume that $0$ is a regular value for $\Phi$
and that $S^1$ acts freely on the level set $\Phi^{-1}(0)$.
Then $Z=\Phi^{-1}(0)$ is a reducible hypersurface, and it
splits $M$ into $M_+ := \Phi^{-1}(0,\infty)$ and
$M_- := \Phi^{-1}(-\infty,0)$.
Conversely, every reducible splitting hypersurface can be obtained
in this way.
\end{Example}


\subsection{Reduction and tangent bundles}

Let $i: Z \to M$ be the inclusion of a reducible hypersurface.
We have the following short exact sequences of vector bundles over $Z$:
\begin{equation} \labell{exact1}
	0 \rightarrow TZ \rightarrow
	i^*(TM) \rightarrow NZ \rightarrow 0
\end{equation}
and
\begin{equation} \labell{exact2}
	0 \rightarrow T\cS \rightarrow TZ \rightarrow
	\pi^*(TM_\red) \rightarrow 0
\end{equation}
where $i$ and $\pi$ are the inclusion and quotient maps,
as in \eqref{reduction}.
Here $T\cS$ is the sub-bundle of $TZ$ consisting of those vectors
that are tangent to the $S^1$ orbits; it is
an oriented trivial real line bundle.
Similarly, $NZ$, the normal bundle to $Z$ in $M$,
is an oriented trivial real line bundle, because $Z$ is co-oriented.

It follows from \eqref{exact1} and \eqref{exact2}
that there exists a (non-canonical) isomorphism of vector bundles over $Z$:
\begin{equation}
\labell{splitting}
   i^* (TM) \cong \pi^* (TM_\red) \oplus T\cS \oplus NZ .
\end{equation}
This isomorphism is $S^1$-equivariant with respect to the natural
$S^1$ action on the bundle $\pi^* (TM_\red)$
and the trivial $S^1$ action on $T\cS \oplus NZ$.
If we choose trivializations of
$NZ$ and $T\cS$, \eqref{splitting} gives an equivariant isomorphism:
\begin{equation} \labell{plus R2}
  i^* (TM) \cong \pi^* (TM_\red) \oplus \R^2
\end{equation}
where $\R^2$ denotes the trivial bundle over $Z$ with fiber $\R^2$.


\subsection{Reduction of orientations}
\labell{red-o}
An orientation $\cal O$ on $M$ determines a {\bf reduced orientation}
${\cal O}_\red$ on $M_\red$ by requiring the induced orientations
on the left and right sides of \eqref{splitting} to agree.

This convention guarantees that if $\cO$ is the orientation
induced by a symplectic form on $M$ and $Z$ is the zero level
set of a moment map, then $\cO_\red$ is the orientation induced by the
reduced symplectic form (see \S\ref{red-twoform}).


\subsection{Reduction of line bundles}
\labell{red-line-bundles}
An equivariant complex line bundle $L$ over $M$ determines a {\bf
reduced line bundle}
$L_\red:= i^*L /S^1$ over $M_\red := Z/S^1$, where $i :Z \to M$ is the
inclusion of a reducible hypersurface.


\subsection{Reduction of almost complex structures}
\labell{red-a-c}
Almost complex structures do not always descend to the reduced
space, but they do when the following condition is satisfied.
\begin{Definition}
\labell{def-admissible}
A reducible hypersurface $Z$ is {\bf admissible} for
an invariant almost complex structure $J$ on $M$
if  $J\xi_M$ is transverse to $Z$ at all points of $Z$,
where $\xi_M$ is the vector field on $M$ that generates the circle action.
\end{Definition}

Given  an admissible hypersurface $Z$ for an invariant almost
complex structure  $J$,
one can construct an almost complex structure on the reduced space
$M_\red = Z / S^1$ as follows.
Let ${\cal V} := TZ \cap J(TZ) \subset i^*(TM)$
be the maximal complex sub-bundle of $TZ$.
The composition
${\cal V} \hookrightarrow TZ \stackrel{\pi_*}{\to} \pi^* (TM_\red)$
is an isomorphism.  This isomorphism defines an equivariant
complex structure on the vector bundle $\pi^* (TM_\red)$,
which in turn descends to an almost complex structure $J_\red$
on $M_\red$.

If $Z$ is an admissible hypersurface, we assign a sign
to each connected component, $X$, of the reduced space, $M_\red:= Z/S^1$:
\begin{equation} \labell{sign}
   \sigma _J(X) = \begin{cases}
	\phantom{-}1 &	\mbox{if $J\xi_M$
			points in the positive direction
			at points of $\pi \inv (X)$}\\
	         -1  &	\mbox{if $J\xi_M$
			points in the negative direction
			at points of $\pi \inv (X)$}
\end{cases}
\end{equation}
where ``positive" and ``negative" are with respect to
the given co-orientation of $Z$.

\begin{Example}
Suppose $(M,\omega)$ is a symplectic manifold
and $J$ is a compatible almost complex structure.
If a level set of the moment map
for a Hamiltonian circle action is reducible,
then it is admissible, and $\sigma_J(M_\red) = +1$.

More generally, let $\omega$ be presymplectic
and let $(M,J,\omega)$ have the property that for every $v \in TM$,
$\omega(v,Jv) = 0$ implies $\iota(v) \omega =0$;
equivalently, ``if $v$ has any friends, $Jv$ is a friend''.
This property holds if,  for example,
$(M,J)$ is the complex blowup of a K\"ahler manifold
along some complex submanifold and $\omega$ is the pullback
of the K\"ahler form.
Again, for any Hamiltonian circle action,
every reducible level set is admissible.

\end{Example}

\begin{Lemma} \labell{orientations}
Let $M$ be a manifold with a circle action,
$J$ an invariant complex structure,
and $\cO$ the induced orientation.
Let $Z$ be an admissible hypersurface, and $X$
a connected component of the reduced space $M_\red = Z / S^1$.
Then the reduced complex structure, $J_\red$,
is compatible with the reduced orientation, $\cO_\red$, on $X$
if and only if $\sigma_J(X) =1$.
\end{Lemma}

\begin{pf}
By the definitions of the reduced complex structure, a positive frame
on $M_\red$ with respect to the complex orientation $J_\red$
can be written as the push-forward of
\begin{equation} \labell{frame}
	\eta_1, J\eta_1, \ldots, \eta_d, J\eta_d
\end{equation}
where $\{\eta_i\}$ is a complex basis  for the maximal complex sub-bundle of
$TZ$.
The frame
$$
	\eta_1, J\eta_1, \ldots, \eta_d, J\eta_d,  \xi_M, J \xi_M
$$
on $M$ is positive with respect to the complex orientation.
By the definition of the reduced orientation,
the push-forward of \eqref{frame}
is positive with respect to $\cO_\red$
if and only if $(\xi_M, J \xi_M) $ is positive in $T\cS \oplus NZ$.
\end{pf}


\subsection{Reduction of stable complex structures}
\labell{red-stable}

Although an almost complex structure often fails to descend to the reduced
space,
it induces a stable complex structure, which does descend to the reduced space.

Let $M$ be a manifold with circle action, and let $Z$ be a reducible
hypersurface.  Take a stable complex structure $\cJ$ represented by
a complex structure on $TM \oplus \R^k$.
By \eqref{plus R2}, this complex structure restricts
to an invariant complex structure on
$\pi^* (TM_\red) \oplus \R^2 \oplus \R^k$,
which in turn descends to a complex structure on $TM_\red \oplus \R^{2+k}$,
providing a stable complex structure ${\cal J}_\red$ on $M_\red$.
The reduced stable complex structure thus obtained will not depend
on either the choice of isomorphism in  \eqref{plus R2}
or the choice of representative of the equivariant stable complex structure.

\begin{Lemma} \labell{a-c-s-c}
Let $S^1$ act on a manifold $M$,
let $Z$ be an admissible hypersurface for an equivariant
almost complex structure $J$, and let
$\cJ = [J]$ denote the corresponding stable complex structure.
Then $\cJ_\red = [J_\red]$.
\end{Lemma}

\begin{pf}
Denote the maximal complex sub-bundle of $TZ$  by $\cal V$.
As complex bundles,
 $i^* TM \cong {\cal V} \oplus \C$,
where the trivial sub-line-bundle $\C$ inside the complex
vector bundle $i^*(TM)$ has the trivializing section $\xi_M$.
Therefore, $({\cal V}, J)$ and $(i^* TM, J)$ are stably equivalent.
The Lemma follows.
\end{pf}


\subsection{Reduction of spin${}^{\bs c}$ structures}
\labell{sec:spinc-reduction}
Let $M$ be an oriented Riemannian manifold with a circle action.
Let $(P,p)$ be an $S^1$-equivariant spin$^c$
structure on $M$ with associated line bundle $L$.
Let $i : Z \to M$ be the inclusion map of a reducible hypersurface
$Z$ into $M$.

By equation \eqref{plus R2}, which reads
$$
	i^* (TM) \cong \pi^* (TM_\red) \oplus \R^2,
$$
the invariant Riemannian metric on $M$
induces an invariant metric on $\pi^* (TM_\red)$,
hence a metric on $M_\red$.
The $\spinc$ structure on $M$ restricts to an $S^1$-equivariant spin$^c$
structure, $(i^*P,i^*p)$, on the vector bundle $i^* TM$.
By \eqref{plus R2} combined with Lemma~\ref{reduced-spin-c},
this induces an $S^1$-equivariant spin$^c$ structure on $\pi^*(TM_\red)$
with the same associated line bundle, $i^*L$.
Taking the quotient by the $S^1$-action yields a spin$^c$ structure on
$M_\red$ with associated line bundle $L_\red$.
This is independent of the choice of isomorphism in \eqref{plus R2}.

\begin{Remark}
The reduction of a $\spinc$ structure induced from a stable complex
structure, a line bundle, an orientation, and a metric is isomorphic
to the $\spinc$ structure induced from the reduced stable complex
structure, reduced line bundle, reduced orientation, and reduced
metric.
\end{Remark}


\subsection{Presymplectic reduction}
\labell{red-twoform}
Let $(M,\omega)$ be a presymplectic manifold with a circle action
and moment map $\Phi : M \to \R$ (see Appendix~\ref{conventions}).

If $a$ is a regular value for $\Phi$, the circle action
on the level set $\Phi\inv(a)$ is locally free. If this action is free,
$\Phi\inv(a)$ is a reducible splitting hypersurface.
There then exists a unique presymplectic
form $\omega^a_\red$ on $M^a_\red := \Phi\inv(a)/{S^1}$
such that $\pi^*\omega_\red = i^* \omega$.

Unlike in the symplectic case, the reduced space $M^a_\red$ need not
be connected even if $M$ is connected and $\Phi$ is proper.

If $[\omega + u \Phi ]$ is the equivariant first Chern class of
$L$ (see \S\ref{presymplectic quantization}),
then $[\omega^a_\red]$ is the first Chern class of $L_\red^a$ if we reduce
at $a = 0$, but not otherwise.


\subsection{Reduction and group actions}
Let $M$ be a manifold with a circle action,
and let $Z$ be a reducible hypersurface.
Suppose that a second group, $G$, acts on $M$,
and that this second action commutes with the circle action
and preserves $Z$. Then we get a $G$-action on the
reduced space, $M_\red$.
For each of the structures discussed in this section,
if we begin with an $(S^1\times G)$-equivariant structure on $M$,
reduction produces a $G$-equivariant structure on $M_\red$.


\section{The main theorems}
\labell{sec:main}

In this section we state three versions of
``quantization commutes with reduction''
corresponding to almost complex, stable complex,
and $\spinc$ quantizations.
In each case, the theorem is written as an equality between
virtual vector spaces, and  holds when the quantization
and reduction data are compatible in the following sense:
if the weight by which $S^1$ acts on a fiber of the line bundle
associated to the quantization data is large (small) over a fixed point,
then that fixed point should be in $M^+$ ($M^-$).
Finally, we prove that the first two of these theorems
follow from the third.

\begin{Theorem}
{\bf Almost complex quantization commutes with reduction,
up to sign.}
\labell{theorem-a-c}
Let the circle act on a compact manifold, $M$.
Let $L$ be an equivariant line bundle over $M$, and
let $J$ be an invariant almost complex structure on $M$.
Let $Z \subset M$ be an admissible reducible hypersurface
that splits $M$ into $M_+$ and $M_-$
(see Definitions~\ref{def-reducible},
\ref{def-splits}, and~\ref{def-admissible}).
Assume that the following conditions are satisfied
for every component $F$ of the fixed-point set:
\begin{equation}
\labell{conditions-a-c}
   \begin{array}{lcl}
   \mbox{ fiber weight at } F 	\geq
   - \sum \mbox{ negative tangent weights at } F
   & \Rightarrow & F \subset M_+ \\
   \mbox{ fiber weight at } F 	\leq
   - \sum \mbox{ positive tangent weights at } F
   & \Rightarrow & F \subset M_-
   \end{array}
\end{equation}
(see Remark \ref{weights1}). Then
\begin{equation}
\labell{formula-a-c}
	Q(M,L,J)^{S^1} \, = \, \sum_{X \subseteq M_\red}
		\sigma_J(X) \cdot Q(X,L_\red,J_\red),
\end{equation}
where the sum is taken over all connected components $X$
of the reduced space $M_\red = Z/{S^1}$,
$\sigma_J(X)$ is the sign associated to $X$ as in \eqref{sign},
$L_\red$ is the reduced line bundle (\S\ref{red-line-bundles}),
$J_\red$ is the reduced almost complex structure (\S\ref{red-a-c}),
and $Q$ denotes almost complex quantization
(\S\ref{almost complex quantization}).
\end{Theorem}

\begin{Remark}
\labell{weights1}
The {\bf fiber weight} at $F$ is the weight by which $S^1$ acts
on a fiber of $L$ over $F$. The {\bf tangent weights} at $F$ are
the weights by which $S^1$ acts on the complex vector space
$(T_mM, J_m)$ for $m$ in $F$.
\end{Remark}

\begin{Remark}
The sign $\sigma_J(X)$ in Theorem~\ref{theorem-a-c}
comes from the discrepancy between the quotient orientation
and the almost complex orientations. See Lemma~\ref{orientations}.
\end{Remark}

\begin{Theorem}
{\bf Stable complex quantization commutes with reduction.}
\labell{theorem-s-c}
Let the circle act on a compact manifold $M$.
Let $L$ be an equivariant line bundle over $M$,
let $\cJ$ be an equivariant stable complex structure on $M$,
and let $\cO$ be an orientation on $M$.
Let $Z$ be a reducible hypersurface that
splits $M$ into $M_+$ and $M_-$
(see Definitions~\ref{def-reducible} and \ref{def-splits}).
Assume that the following conditions are satisfied
for every component $F$ of the fixed-point set:
\begin{equation}
\labell{conditions-s-c}
   \begin{array}{lcl}
   \mbox{ fiber weight at } F 	\geq
   - \sum \mbox{ negative tangent weights at } F
   & \Rightarrow & F \subset M_+  \\
   \mbox{ fiber weight at }  F	\leq
   - \sum \mbox{ positive tangent weights at }F
   & \Rightarrow & F \subset M_-
   \end{array}
\end{equation}
(see Remark \ref{weights2}).  Then
\begin{equation}
\labell{formula-s-c}
	Q(M,L,\cJ,\cO)^{S^1} \, = \,
		Q(M_\red,L_\red,\cJ_\red,\cO_\red),
\end{equation}
where $M_\red = Z/S^1$ is the reduced space,
$L_\red$ is the reduced line bundle (\S\ref{red-line-bundles}),
$\cJ_\red$ is the reduced stable complex structure (\S\ref{red-stable}),
$\cO_\red$ is the reduced orientation (\S\ref{red-o}), and
$Q$ denotes stable complex quantization
(\S\ref{stable complex quantization}).
\end{Theorem}

\begin{Remark}
\labell{weights2}
Again, the fiber weight at $F$ means the weight by which
$S^1$ acts on a fiber of $L$ over $F$.
The negative and positive tangent weights of a stable complex
$\cJ$ which is represented  by a complex vector bundle $(TM \oplus \R^k, J)$,
are the negative and positive weights by which $S^1$ acts
on the complex vector space $(T_mM \oplus \R^k, J_m)$ for $m \in F$.
These weights are well defined, even though different representatives
of the stable complex structure may possess different  numbers of zero
weights.
\end{Remark}

We now state our ``main theorem'',
which implies Theorems 1 and 2.
It will be proved in section \ref{proof-spinc}.

\begin{Theorem}
\labell{theorem-spinc}
{\bf Main theorem: $\bs \spinc$ quantization commutes with reduction.}
Let the circle act on a compact oriented Riemannian manifold $M$.
Let $(P,p)$ be an equivariant spin$^c$ structure on $M$.
Let $Z$ be a reducible hypersurface that
splits $M$ into $M_+$ and $M_-$
(see Definitions~\ref{def-reducible} and \ref{def-splits}).
Assume that the following conditions are satisfied
for every component $F$ of the fixed-point set:
\begin{equation} \labell{conditions-spinc}
   \begin{array}{lcl}
   \mbox{ fiber weight at } F 	\geq
   \phantom{-} \sum | \mbox{tangent weights at }F|
   & \Rightarrow & F \subset M_+ \\
   \mbox{ fiber weight at } F 	\leq
   - \sum | \mbox{tangent weights at } F|
   & \Rightarrow & F \subset M_-
   \end{array}
\end{equation}
(see Remark \ref{weights3}). Then
\begin{equation}
\labell{formula-spinc}
	Q(M,P,p)^{S^1} \, = \,
		Q(M_\red,P_\red,p_\red),
\end{equation}
where $M_\red = Z/{S^1}$ is the reduced space,
$(P_\red,p_\red)$ is the reduced spin$^c$ structure
(\S\ref{sec:spinc-reduction}),
and $Q$ denotes spin$^c$ quantization (\S\ref{spinc quantization}).
\end{Theorem}

\begin{Remark}\labell{weights3}
The fiber weight at $F$ means the weight by which $S^1$
acts on the fiber over $F$ of the line bundle associated
to the $\spinc$ structure.
The  tangent weights are defined only up to sign, because
the vector space $T_mM$, for $m\in F$, is only a {\em real\/}
representation space for $S^1$.
However, their absolute values are well defined.
\end{Remark}


\begin{pf*}{Proof that Theorem~\ref{theorem-a-c} follows from
Theorem~\ref{theorem-s-c} }
We will now show that the almost complex version of
``quantization commutes with reduction'' follows from
the stable complex version.

Let $M$, $L$, $J$, and $Z$ be as in the statement
of Theorem~\ref{theorem-a-c}.
Let $\cJ$ be the equivariant stable complex structure induced by $J$,
and let $\cO$ be the orientation induced by $J$.
Conditions \eqref{conditions-a-c} and \eqref{conditions-s-c}
are the same, because the nonzero tangent weights for $\cJ$
are the same as those for $J$.
The corresponding quantizations are also the same:
$Q(M,L,J) = Q(M,L,\cJ,\cO)$, and
the almost complex quantization of the reduced
space is equal to the stable complex quantization
of the stable complex structure and orientation
induced by $J_\red$ (see Remark~\ref{a-c-q = spinc-q}).
The stable complex structure induced by $J_\red$
coincides with $\cJ_\red$ by Lemma~\ref{a-c-s-c}.
However, on a component $X$ of $M_\red$,
the complex orientation of $J_\red$
coincides with the reduced orientation $\cO_\red$
if and only if $\sigma_J(X) = + 1$ (see Lemma~\ref{orientations}).
Since switching orientation results in switching the sign
of the quantization,
$Q(X,L_\red,\cJ_\red,\cO_\red) = \sigma_J(X) Q(X,L_\red,J_\red).$

\end{pf*}


\begin{pf*}{Proof that Theorem~\ref{theorem-s-c} follows from
Theorem~\ref{theorem-spinc}}
We will now show that
the stable complex version of ``quantization commutes with reduction''
follows from the  $\spinc$ version.

The data $(M,L,\cJ,\cO)$ (together with an invariant metric)
determines an equivariant spin$^c$ structure $(P,p)$ on $M$
(see \S\ref{stable complex quantization}).

First, we will show that Conditions \eqref{conditions-s-c}
for the stable complex structure imply Conditions
\eqref{conditions-spinc} for the associated spin$^c$ structure.
The line bundle associated to the spin$^c$ structure
is $ L \otimes L \otimes K^* $, where $K^*$ is the dual to the
canonical line bundle of the stable complex structure.
Thus, over any component of the fixed point set,
\begin{eqnarray*}
	\mbox{the fiber weight of $L \otimes L \otimes K^* $ }
	 & = & 2 \times (\mbox{the fiber weight of $L$}) \\
 	 &   & + \sum \mbox{the tangent weights} \ ,
\end{eqnarray*}
and
\begin{eqnarray*}
	\sum |\mbox{the tangent weights}|
	& = & - 2 \times \sum \mbox{the negative tangent weights}  \\
	& & + \sum \mbox{the tangent weights}\ .
\end{eqnarray*}
Hence the first condition in \eqref{conditions-spinc},
which reads
$$
	\mbox{the fiber weight of $L \otimes L \otimes K^*$}
	\geq \sum |\mbox{the tangent weights}| \ ,
$$
is equivalent to the first condition in \eqref{conditions-s-c},
which reads
$$
	\mbox{the fiber weight of $L$}
	\geq - \sum \mbox{the negative tangent weights}\ .
$$
A similar argument shows that the second conditions
in \eqref{conditions-s-c} and \eqref{conditions-spinc}
are equivalent.

Moreover, by the definition of the quantization  of
a stable complex structure,
the corresponding equivariant quantizations are the same,
$Q(M,L,\cJ,\cO) = Q(M,P,p)$, and the quantization of
the $\spinc$ structure induced by $\cJ_\red$ and
$\cO_\red$ is equal to the stable complex
quantization of $\cJ_\red$ and $\cO_\red$.
Finally, the spin$^c$ structure associated to
$(M_\red,L_\red,\cJ_\red,\cO_\red)$ (and the reduced metric)
is isomorphic to the  reduced spin$^c$ structure, $(P_\red,p_\red)$.
\end{pf*}

\begin{Remark}
\labell{vary0}
The criteria  which the hypersurface must satisfy at first seem surprisingly
lax.   However, this is much less surprising when  one considers the following
facts.
Assume, for the sake of avoiding orbifolds, that the action of $S^1$
on $M$ is semi-free, that is,  every point is fixed or free.
Fix some quantizing information, and
consider any splitting hypersurface $Z$.
Let $\{F_i\}_{i=1}^l$ be the set of components of the fixed-point set
contained in $M^+$.
Let $B_i$ be a small tubular neighborhood of $F_i$, and let $S_i$ be its
boundary.
Then $Z/{S^1}$  is cobordant to $\coprod_{i=1}^l S_i/{S^1}$.
(This cobordism was used by Victor Guillemin and by Shaun Martin.)
Moreover, the cobording manifold, $(M \ssminus \sqcup_i B_i) / S^1$,
carries the quantization data.
Since  quantization is a cobordism invariant, this shows that
$Q(Z/{S^1}) = \sum_{i=1}^l Q(S_i/{S^1})$, and hence explains part
of the freedom.  The rest is explained once we notice that if
$F_i$ is a component of the fixed-point set which is allowed to be
in either $M^+$ or $M^-$, the quantization of  $S_i/{S^1}$  is empty.
\end{Remark}

\begin{Example}
The four-sphere, $S^4$,
has no symplectic or almost complex structure (see \cite{BS,karoubi}).
However, the inclusion
$$
   j: S^4\ \hookrightarrow \C^3\ , \ \ \
   j(S^4) = \{ (z_1,z_2,z_3) : |z_1|^2+|z_2|^2+(\Re z_3)^2 = 1
                               \mbox{ and } \Im z_3 = 0 \},
$$
induces a stable complex structure,
$$
   TS^4 \oplus \R^2 \cong j^*(T\C^3) \cong \C^3,
$$
where the $\R^2$ is the trivial normal bundle of $S^4$
inside $\C^3$.
Equip $S^4$ with the standard orientation as a subset of
$\R^5 \cong \{ (z_1,z_2,z_3) : \Im z_3 = 0 \}$,

Let $L = \C$ be the trivial line bundle over $S^4$ equipped with
the following $S^1$ action:
$$
   e^{i\theta} \cdot \left( (z_1,z_2,z_3),v \right)
   = \left( (e^{i\theta}z_1,e^{i\theta}z_2,z_3),e^{mi\theta}v \right)
$$
for $(z_1,z_2,z_3) \in S^4$ and $v \in \C$.
The equator, $Z = \{ (z_1,z_2,0) : |z_1|^2+|z_2|^2= 1 \} \cong S^3$,
is a reducible hypersurface that splits $M$ into the northern and
southern hemispheres.
At each of the two fixed points, $N=(0,0,1)$ and $S=(0,0,-1)$,
the non-zero weights of the stable tangent bundle are $+1$ and $+1$,
and the fiber weight of $L$ is $m$.
Therefore, the conditions of  Theorem~\ref{theorem-s-c}
are satisfied exactly when $m =-1$.

We can now check by computing explicitly:
The orientation on $S^4$ agrees with that of
$\C^2 \cong \{ (z_1,z_2,0) \in \C^3 \}$ at $N$, but is opposite at $S$.
Therefore, the two fixed-point contributions
to the index-character formula cancel each other.
Hence,  $\dim Q(M)^{S^1} = 0$ for all $m$.
On the other hand,  $M_\red = S^2$,
$\cJ_\red$ is the class of the standard complex structure,
$\cO_\red$ is also standard, and $L_\red = \cO (m)$.
Consequently, $\dim Q(M_\red) = m+1$.
As required, $Q(M)^{S^1} = Q(M_\red)$ when $m = -1$.
\end{Example}

\begin{Example}
Let $S^1$ act on any spin manifold $M$, and  consider the $\spinc$
bundle with trivial associated line bundle.  The criteria are
automatically satisfied on every fixed point set, so quantization
commutes with reduction for every splitting hypersurface.
In this case, both are trivial.
\end{Example}

\begin{Example}
Notice that our theorem makes sense whenever $Z$ is a reducible hypersurface.
However, even if $M \ssminus Z$ can be written as the disjoint union of
an $M^-$ and an $M^+$ which satisfy conditions~\eqref{conditions-a-c},
{}~\eqref{conditions-s-c}, or ~\eqref{conditions-spinc},
the theorem is generally not true unless $Z$ is splitting:

Consider the torus $T = \R^2/\Z^2$ with the complex structure induced
by the isomorphism $\R^2 = \C$ and the trivial line bundle.
The circle $\R/\Z$ acts on $T$ by  $[\alpha] \cdot [x,y] = [x + \alpha, y]$.

Because there are no fixed points, by Atiyah-Segal-Singer the almost
complex quantization of $T$ is the zero representation.
In particular, the zero weight space is zero.

On the other hand, the circle $Z := [x,0]$ is a reducible hypersurface.
If we take $M^+ := T \ssminus  Z$ and $M^- = \emptyset$,
then they satisfy conditions~\ref{conditions-a-c} since neither
contains any fixed point.
Nevertheless, the reduced space is a point, so its quantization is not zero.
Quantization does not commute with reduction in this example,
because  $Z$ is not splitting.
\end{Example}


\section{Presymplectic corollaries}
\labell{sec:corollaries}

The most important application of our theorems is to
presymplectic manifolds, where  we quantize and reduce with
respect to the same closed two form and moment map.
In this case, the compatibility conditions in the earlier theorems
are automatically satisfied, and the appropriate version of
``quantization commutes with reduction'' holds at every regular integer value.

More precisely, let $S^1$ act on a manifold $M$.
Let $\omega$ be a presymplectic form, that is, a closed invariant two-form
with a moment map $\Phi : M \to \R$.
Assume that $\omega + u \Phi$ is prequantizable,
{\em i.e.\/}, that either $\omega + u\Phi$ is complex prequantizable
and that $M$ has an invariant almost complex structure or an
invariant stable complex structure and an orientation,
or that $M$ is oriented and $\omega + u \Phi$ is $\spinc$ prequantizable
(see \S\ref{presymplectic quantization}).
If zero is a regular value of $\Phi$,
then $Z = \Phi\inv(0)$ is automatically a reducible splitting hypersurface.
Moreover, the compatibility conditions --
Conditions \eqref{conditions-a-c}, \eqref{conditions-s-c}, or
\eqref{conditions-spinc} -- are automatically satisfied:
the fiber weight at $F$ is equal to either $\Phi(F)$ (in the almost complex
and stable complex cases) or $2\Phi(F)$ (in the $\spinc$ case);
this is positive if $F$ is in $M_+$
and negative if $F$ is in $M_-$.
The compatibility conditions follow immediately
by comparing the signs of the left and right terms
in the inequalities
\eqref{conditions-a-c}, \eqref{conditions-s-c}, or \eqref{conditions-spinc}.
This shows that quantization commutes with reduction at zero.

More generally, consider any regular $a \in \Z$.
Then the reduction $\omega_\red^a$ of  $[\omega + u \Phi]$ at $a \in \R$
is identical to the reduction of $[\omega + u( \Phi - a)]$ at $0$.
Therefore, by the preceding paragraph, the multiplicity of $0$ in
$Q(M,\omega, \Phi-a, \ldots )$ is
$\dim \ Q(M_\red^a,\omega_\red^a,\ldots)$.
But, by Remark~\ref{shift}, the multiplicity of $0$ in
$Q(M,\omega,\Phi-a,\ldots)$ is simply the multiplicity of $a$
in $Q(M,\omega,\Phi,\ldots)$.
Theorems \ref{theorem-a-c}--\ref{theorem-spinc} thus have the
following corollaries.

\begin{Corollary}
{\bf Almost complex quantization of presymplectic manifolds
commutes with reduction, up to sign.}
\labell{coroll1}
Let $S^1$ act on a compact presymplectic and almost complex
manifold $(M,\omega, J)$
with a moment map $\Phi$. Assume that $\omega$ and $\Phi$ are integral
(see \S\ref{presymplectic quantization}).
Let $a \in \Z$ be a regular value for $\Phi$ such that the $S^1$
action on the $a$-level set is free and  the $a$-level set
is admissible (see Definition \ref{def-admissible}). Then
\begin{equation} \labell{eq:cor1}
	Q(M,\omega,\Phi,J)^{a} \, = \, \sum_{X \subseteq M^a_\red}
		\sigma_J(X) \cdot Q(X,\omega^a_\red,J^a_\red),
\end{equation}
where we sum over the connected components $X$
of the reduced space $M^a_\red := \Phi^{-1}(a)/{S^1}$,
$\sigma_J(X) = \pm 1$ is the sign associated to $X$ as in \eqref{sign},
$\omega^a_\red$ is the reduced presymplectic form (\S\ref{red-twoform}),
$J^a_\red$ is the reduced almost complex structure (\S\ref{red-a-c}),
$Q$ denotes the almost complex quantization (\S\ref{presymplectic
quantization}),
and $Q(\cdot)^a$ is the $a$th weight space in $Q(\cdot)$.
\end{Corollary}

\begin{Remark}
When $\omega$ is symplectic and $J$ is a compatible almost complex structure,
every reducible hypersurface is admissible, and $\sigma_J(M_\red^a) =1$.
In this case, Corollary~\ref{coroll1} reduces to the usual theorem,
that  ``quantization commutes with reduction" for
almost K\"ahler quantization on symplectic manifolds.
\end{Remark}

\begin{Corollary}
\labell{coroll2}
{\bf Stable complex quantization of presymplectic manifolds commutes
with reduction.}
Let $S^1$ act on a compact oriented presymplectic and stable
complex manifold $(M,\omega,\cJ,\cO)$
with a moment map $\Phi$.  Assume that $\omega$ and $\Phi$ are integral
(see \S\ref{presymplectic quantization}).
Assume that $a \in \Z$ is a regular value for $\Phi$ and that the $S^1$ action
on the $a$-level set is free.
Let $\cJ$ be an equivariant stable complex structure.  Then
\begin{equation} \labell{eq:cor2}
	Q(M,\omega,\Phi,\cJ,\cO)^{a} \, = \,
		Q(M^a_\red,\omega^a_\red,\cJ^a_\red,\cO^a_\red),
\end{equation}
where $M^a_\red := \Phi^{-1}(a)/{S^1}$ is the reduced space,
$\omega^a_\red$ is the reduced presymplectic form (\S\ref{red-twoform}),
$\cO^a_\red$ is the reduced orientation (\S\ref{red-o}),
$\cJ^a_\red$ is the reduced stable complex structure (\S\ref{red-stable}),
$Q$ denotes the stable complex quantization
(\S\ref{stable complex quantization} and \S\ref{presymplectic quantization}),
and $Q(\cdot)^a$ is the $a$th weight space in $Q(\cdot)$.
\end{Corollary}

\begin{Corollary}
\labell{coroll3}
{\bf Spin$^c$ quantization of presymplectic manifolds commutes with
reduction.}
Let $S^1$ act on a compact presymplectic manifold, $(M,\omega)$,
with a moment map $\Phi$. Assume that $(M,\omega,\Phi)$ is equivariantly
spin$^c$ prequantizable (see \S~\ref{presymplectic quantization}).
Assume that $a \in \Z$ is a regular value for $\Phi$ and that the $S^1$ action
on the $a$-level set is free. Let $\cO$ be an orientation on $M$.  Then
\begin{equation} \labell{eq:cor3}
	Q(M,\omega,\Phi,\cO)^{S^1} \, = \,
		Q(M^a_\red,\omega^a_\red,\cO^a_\red),
\end{equation}
where  $M^a_\red = \Phi\inv(a) / S^1$ is the reduced space,
$\omega^a_\red$ is the reduced presymplectic form (\S\ref{red-twoform}),
$\cO^a_\red$ is the reduced orientation (\S\ref{red-o}),
$Q$ denotes the spin$^c$ quantization
(\S\ref{spinc quantization} and \S\ref{presymplectic quantization}),
and $Q(\cdot)^a$ is the $a$th weight space in $Q(\cdot)$.
\end{Corollary}

\begin{Remark} \labell{vary1}
Assume that  $0$ is a regular value of $\Phi$.
Since tangent weights are always integers,
$t$ is also a regular value for all $-1 < t < 1$
and conditions~\ref{conditions-a-c}  also hold for $\Phi\inv(t)$.
Thus, Formula \eqref{formula-a-c} holds for
$M_\red^t = \Phi\inv(t) / S^1$.
(This does not follow from Corollary~\ref{coroll1}.)
A similar result holds in the stable complex case.

A much stronger, and more useful, result holds in the $\spinc$ case:
even if $0$ is a singular value, conditions \eqref{conditions-spinc} hold
for $Z = \Phi\inv(t)$ for all $ -\half < t < \half$.
(The boundaries are $\pm \half$ and not $\pm 1$
because the fiber weight at $F$ is $2\Phi(F)$, not $\Phi(F)$.
See \S\ref{presymplectic quantization}.)
Thus Formula \eqref{formula-spinc} holds for any
$M_\red^t = \Phi\inv(t) / S^1$ with $ -\half < t < \half$.
The manifold $M_\red^t$ then provides a ``resolution''
of $\Phi\inv(0) / S^1$, if the latter is singular, in the spirit
of Meinrenken and Sjamaar \cite{MS}.
Again, this does not follow from Corollary~\ref{coroll3}.
\end{Remark}

\begin{Remark}
Let $S^1$ act on an almost complex manifold $(M,J)$.
Let $L$ be a complex Hermitian line bundle, and let
$K^*$
be the dual of the canonical line bundle.
Let $\omega$ and $\alpha$ be two-forms representing the first Chern class
of $L$ and $K^*$, respectively.
Let $\Phi$ and $\Psi$ be their respective moment maps.
Consider the moment map
$$
	\Phi_t = \Phi + t \Psi
$$
for the presymplectic form
$$
	\omega_t = \omega + t \alpha\ .
$$
Then for all $t \in [0,1]$, if $\Phi_t \inv (0)$ is admissible,
Formula~\eqref{formula-a-c} still holds for
$$
	M_\red^t = \Phi_t \inv (0) / S^1\ .
$$
This follows from Theorem~\ref{theorem-a-c}. The compatibility
conditions, \eqref{conditions-a-c}, are automatically satisfied:
\begin{eqnarray*}
F \subset M_+  & \Leftrightarrow & \Phi_t(F) >0 \\
 & \Leftrightarrow &
\mbox{fiber weight at $F$} + t \sum \mbox{tangent weights at $F$} >0 \\
 & \Rightarrow & \mbox{fiber weight at $F$} > -t \sum \mbox{tangent weights
at $F$} \\
 & & \geq -t \sum \mbox{positive tangent weights at $F$} \\
 & & \geq - \sum \mbox{positive tangent weights at $F$}
\end{eqnarray*}
and, similarly,
$$ F \subset M_- \Rightarrow
\mbox{fiber weight at $F$} < - \sum \mbox{negative tangent weights at $F$}.
$$
Corollary~\ref{coroll1} is the case $t=0$, and
Corollary~\ref{coroll3} (for the associated $\spinc$ structure)
corresponds to the middle point, $t = {1 \over 2}$.
Hence, Corollary~\ref{coroll3} does {\em not\/}
imply Corollary~\ref{coroll1}; they are different theorems.
A similar remark holds for stable complex quantization.
\end{Remark}

\begin{Example}
Let $S^1$ rotate $S^2$, fixing the north and south poles.
Consider the two-form $m \sigma$,
where $m$ is an integer, and $\sigma$ is the standard
K\"ahler form on $S^2 \simeq \CP^1$ with total integral equal to one.
Choose the moment map $\Phi: S^2 \to \R$
so that it takes values  $0$ and $m$ at the south and north poles,
respectively. Equip $S^2$ with the standard complex structure, $J$.

The line bundle with equivariant first Chern class  $[\omega + u \Phi]$
is $\cO(m)$ with fiber weights $0$ and $m$ over the south and north poles,
respectively.

If $0 < m$, the multiplicity of $a$ in the quantization,
$Q(S^2, m\sigma, \Phi,J)$, is $1$ for all integers $a$ such that
$0 \leq a \leq m$, and is $0$ otherwise.
The singular values are $a=0,m$.
For $ 0 < a < m$,  $M^a_\red$ is a single point
with reduced orientation $+1$.
If $a<0$ or $a>m$, $M^a_\red$ is empty.
Hence for any regular integer $a$,
$Q(S^2)^a = Q(M_\red^a)$; the dimension of these spaces is $+1$
if $0 < a < m$ and $0$ if $a < 0$ or $a > m$.

If $m < 0$, the multiplicity of $a$ in the quantization
$Q(S^2, m\sigma, \Phi,J)$ is $-1$ for all integers $a$ such that
$m < a < 0$, and is $0$ otherwise.
The singular values are $a=0,m$.
For $ m < a < 0$,  $M^a_\red$ is a single point
with reduced orientation $-1$.
If $a < m$ or $a > 0$, $M^a_\red$ is empty.
Hence for any regular integer $a$,
$Q(S^2)^a = -Q(M_\red^a)$;
the dimension of these spaces is $-1$ if $m < a < 0$
and $0$ if $a < m$ or $a > 0$.
\end{Example}

\begin{Example}
Let $S^1$ rotate $S^2$ and consider the two-form $m \sigma$, as  in
the previous example.
Choose the moment map $\Phi: S^2 \to \R$
so that it takes values  $-1/2$ and $m - 1/2$ at the south and north poles,
respectively.

If $0 < m$, the multiplicity of $a$ in the spin$^c$
quantization, $Q(S^2,m\sigma, \Phi,\cO)$, is $+1$ for
all integers $a$ such that
$0 \leq a < m$ and $0$ otherwise.
The singular values are $a=-1/2,m-1/2.$
For $-1/2 < a < m-1/2$,  $M_\red$ is a single point
with reduced orientation $+1$.
For $a <-1/2$ or $a > m-1/2$, $M_\red$ is empty.
Hence $Q(S^2)^a = Q(M_\red^a)$ for all integers $a$;
the dimensions of these spaces are $+1$ if
$0 \leq a <m$ and $0$ otherwise.

If $m < 0$, then the multiplicity of $a$ in the quantization
$Q(S^2,m\sigma,\Phi,\cO)$ is $-1$ for all integers $a$ such that
$m \leq a < 0$, and $0$ otherwise.
The singular values are $a=-1/2$ and $a=m-1/2$.
For $ m-1/2 < a < -1/2$,  $M^a_\red$ is a single point,
with reduced orientation $-1$. For $a < m-1/2$ or $a > -1/2$,
$M^a_\red$ is empty. Hence $Q(S^2)^a = Q(M_\red^a)$ for all integers $a$;
the dimensions of these spaces are $-1$ if
$m \leq a <0$ and $0$ otherwise.
\end{Example}

\begin{Example}
Now consider the presymplectic form
$1/2\pi \left( 10z + 5 \right) dz \, d\theta$ with moment map
$\Phi (z,\theta) = 5 z^2 + 5 z + 1$,
where $(z,\theta)$ are cylindrical coordinates on $S^2$,
the unit sphere in $\R^3$.

We will work out the almost complex and the stable complex quantizations.
The quantizing line bundle is $O(10)$, with fiber weight $11$
at the north pole and $1$ at the south pole.
The multiplicity of the
weight $a$ is $+1$ if $1 \leq a \leq 11$
and $0$ otherwise.

The zero-level set has two components,
hence so does the reduced space: $M_\red = X_1 \sqcup X_2$.
Almost complex quantization gives $\dim Q(X_1) = \dim Q(X_2) = +1$.
However, $\sigma_J(X_1) = +1$ and $\sigma_J(X_2) = -1$,
so Corollary \ref{coroll1} gives
$\dim Q(S^2)^{S^1} = \dim Q(X_1) - \dim Q(X_2) = 0 $.

Stable complex quantization for the stable complex structures
associated to these almost complex structures is very similar.
The only difference lies in the definition of the reduced orientation.
It gives $\dim Q(S^2)^{S^1} = \dim Q(X_1) + \dim Q(X_2)$,
with $\dim Q(X_1) = +1$ and $\dim Q(X_2) = -1$.
\end{Example}


\section{The cut space}
\labell{sec:cutting}

Following \cite{4},
to prove that quantization commutes with reduction we
construct a new space by  cutting.
The technique of symplectic cutting
was introduced by Lerman in \cite{L}
for symplectic manifolds with a circle action and a moment map.
We apply Lerman's construction to manifolds with other
structures.
Our goal in this section is to prove the following Proposition.

\begin{Proposition} \labell{cutting}
Let the circle group, $S^1$, act on a smooth manifold $M$.
Let $Z$ be a reducible hypersurface (see Definitions \ref{def-reducible}
and \ref{def-splits}) that splits $M$ into $M_+$ and $M_-$,
and let $M_\red := Z / S^1$ be the reduced space.
\begin{enumerate}
\item
There exists a smooth manifold, $M_\cut$, with a circle action,
and equivariant embeddings
$$
        i_+ : M_+ \hookrightarrow M_\cut
        \quad \mbox{and} \quad
        i_\red: M_\red \hookrightarrow M_\cut \ ,
$$
such that $M_\cut$ is the disjoint union of $i_+ (M_+)$ and $i_\red(M_\red)$.
Here the circle action on $M_+$ is induced from $M$,
and the action on $M_\red$ is trivial.
\item
The normal bundle to the submanifold $i_\red(M_\red)$
in $M_\cut$, pulled back to $M_\red$ via the map $i_\red$,
is isomorphic to the bundle
\begin{equation} \labell{normal}
  N := Z \times_{S^1} \C
\end{equation}
over $M_\red :=Z / S^1$.
Here $[m,z] = [a \cdot m , a \inv z]$ for all $a \in S^1$,
and the circle acts on this bundle, $N$, by $ \lambda \cdot [m,z] = [ m
,\lambda z]$.
\item
Given an orientation on $M$,
there exists an orientation on $M_\cut$ such that

\begin{enumerate}
\item
the open embedding $i_+ : M_+ \hookrightarrow M_\cut$ preserves orientation;
\item
the orientation of the submanifold $i_\red(M_\red)$
coming from the reduced orientation of $M_\red$ (see \S\ref{red-o}),
followed by the orientation on the normal bundle of $i_\red(M_\red)$
that is induced from the complex orientation on the fibers
of \eqref{normal}, gives the orientation on $M_\cut$.
\end{enumerate}

\item
Let $(P,p)$ be an equivariant spin$^c$ structure over $M$
with an associated line bundle $L$.
Then there exists an equivariant spin$^c$ structure $(P_\cut,p_\cut)$ on
$M_\cut$
whose associated line bundle, $L_\cut$, pulls back to the following
$S^1$-equivariant complex line bundles on $M_+$ and on $M_\red$:
\begin{enumerate}
\item
$ i_+^* L_\cut \cong L|_{M_+} $
\item
$ i_\red^* L_\cut \cong L_\red \otimes N $ \quad (see \S\ref{red-line-bundles}
and \eqref{normal}).
\end{enumerate}
\end{enumerate}
\end{Proposition}

\begin{pf*}{Proof of part 1}
We can assume that $Z$ is the zero level set of an invariant function
$\Phi : M \to \R$ for which $0$ is a regular value, and that
$M_+$ is the pre-image under $\Phi$ of the positive numbers.
Consider the product $M \times \C$ with the anti-diagonal
circle action, $a \cdot (m,z) = (a \cdot m, a^{-1} z)$,
and with the function $\tPhi(m,z) := |z|^2 - \Phi(m) $.
It is easy to check that $0$ is a regular value for $\tPhi$
and that $S^1$ acts freely on the zero level set of $\tPhi$,
$$\tZ = \{ (m,z) \ | \ \Phi(m) = |z|^2 \}.$$
The {\bf cut space}, defined by
$$
	M_\cut := \tZ / S^1 = (M \times \C)_\red \ ,
$$
is therefore a smooth manifold.
Let $S^1$ act on $M_\cut$ by
$\lambda \cdot [m,z] = [\lambda \cdot m , z] = [m , \lambda z]$.

The inclusion map $i_+ : M_+ \hookrightarrow M_\cut$ is induced
by the following commuting diagram
\begin{equation} \labell{i-plus}
\begin{array}{ccc}
   M_+ & \stackrel{\tilde{i}_+}{\to} & \tZ \\
   & \stackrel{i_+}{\searrow}
   & \phantom{\scriptstyle \pi_\cut} \downarrow
			{\scriptstyle \pi_\cut} \\
   & & \ M_\cut
\end{array}
\end{equation}
from the map $\tilde{i}_+ (m) := \big(m,\sqrt{\Phi(m)} \big)$.
The inclusion map $i_\red : M_\red \hookrightarrow M_\cut$
is induced by the following commuting diagram
\begin{equation} \labell{i-red}
\begin{array}{ccc}
   Z & \stackrel{\tilde{i}_\red}{\to} & \tZ \\
   \phantom{\scriptstyle \pi} \downarrow {\scriptstyle \pi} & &
   \phantom{\scriptstyle \pi_\cut} \downarrow {\scriptstyle \pi_\cut} \\
   M_\red & \stackrel{i_\red}{\to} & M_\cut
\end{array}
\end{equation}
from the map $\tilde{i}_\red(m) := (m,0)$.
\end{pf*}

\begin{pf*}{Proof of part 2}
The submanifold $i_\red(M_\red)$ of $M_\cut$ is the quotient of
the submanifold $\ti_\red (Z) = Z \times \{0\}$ of $\tZ$.
The normal bundle of the latter is $Z \times \C$,
with the anti-diagonal circle action induced from $M \times \C$.
\end{pf*}

\begin{pf*}{Proof of part 3(a)}
The orientations on $M$ and $\C$ induce an orientation on $M \times \C$
which, by \S\ref{red-o}, descends to an orientation on
$M_\cut = (M \times \C)_\red$.
Recall, a frame on $M_\cut$ is oriented
if and only if its lifting to $\tZ$,
followed by the generator of the anti-diagonal circle action on $\tZ$,
and further followed by a positive normal to $\tZ$,
gives an oriented frame on $M \times \C$.

Let $\eta_1, \ldots, \eta_d$ be an oriented frame on $M_+$.
Its push-forward to $\tZ$ has the form
$$
	(\tiplus)_* \eta_i = \eta_i + a_i \dd{r} \
$$
for some real numbers $a_i$.
(On the right we identify the vectors
$\eta_i \in T_mM$ and $\dd{r} \in T_z\C$
with their images in
$T_{(m,z)} (M \times\C) \simeq T_mM \times T_z\C$.)
The map $i_+$ preserves orientation if and only if the following
frame on $M \times \C$ is oriented:
\begin{equation} \labell{basis}
	 \eta_1 + a_1 \dd{r}, \ldots, \eta_d + a_d \dd{r},
	 \xi_M - \dd{\theta},\dd{r}
\end{equation}
where $r$, $\theta$ are polar coordinates on $\C$;
this is because $\xi_M - \dd{\theta}$ generates the anti-diagonal circle
action,
and $\dd{r}$ is a positive normal to $\tZ$.
Finally, it is easy to check that \eqref{basis} is an oriented
frame on $M \times \C$.
\end{pf*}

\begin{pf*}{Proof of part 3(b)}
Fix a point $m$ in $Z$.
Consider an oriented basis for $T_mM$ of the form
\begin{equation} \labell{basis1}
	\zeta_1, \ldots, \zeta_d,  \xi_M, n
\end{equation}
where $\zeta_i$ are tangent to $T_mZ$,
where $\xi_M$ generates the circle action,
and where $n$ is a positive normal vector to $Z$ in $M$.
By definition of the reduced orientation, the vectors
$ \zeta_1, \ldots, \zeta_d $
descend to an oriented frame on $M_\red$.

The vectors
\begin{equation} \labell{basis2}
	\zeta_1, \ldots, \zeta_d, \dd{x}, \dd{y}, \xi_M ,n
\end{equation}
form an oriented basis of $T_{(m,0)} (M \times \C)$,
where $x$, $y$ are the Euclidean coordinates on $\C$.
(Here we identified the vectors $\zeta_i, \xi_M, n \in T_mM$
and $\dd{x}, \dd{y} \in T_0\C$ with their images in
$T_{(m,0)} (M \times \C) = T_mM \times T_0 \C$.)
Again by definition of the reduced orientation, the vectors
$$
	\zeta_1, \ldots, \zeta_d, \dd{x}, \dd{y}
$$
descend to an oriented frame on $M_\cut$.
Since the first $d$ of these descend to an oriented frame on
the submanifold $i_\red(M_\red)$, and the last two descend
to an oriented basis of the fiber of the normal bundle
when identified with $\C$,
we conclude that the orientation of the submanifold $i_\red(M_\red)$
followed by that of its normal bundle gives the orientation of $M_\cut$.
\end{pf*}

\begin{pf*}{Proof of part 4}
We will now describe a $\spinc$ structure on $M_\cut$.
First, noting that $\spinc(2) = U(1) \times_{\Z_2} U(1)$,
we define an equivariant $\spinc$ structure over $\C$
by the trivial principal bundle
$$
	\C \times (U(1) \times_{\Z_2} U(1)) \to \C
$$
with the non-trivial circle action
$$ e^{i\theta} \cdot (z,[a,b])
   = (e^{i\theta} z , [e^{i \theta /2} a ,
      		e^{ - i \theta /2 } b]).$$
The associated line bundle is the trivial line bundle $C$
over $\C$ with the circle action
\begin{equation} \labell{canonical}
  e^{i\theta} \cdot (z,b)
    = (e^{i\theta} z, e^{-i\theta} b ).
\end{equation}
In Lemma~\ref{prodspinc} below we define a spin$^c$
structure on the product of any two spin$^c$ manifolds.
In this way we get a $\spinc$ structure  on $M \times \C$,
which in turn descends to a $\spinc$ structure on
$M_\cut = ( M \times \C ) _\red$ (see \S\ref{sec:spinc-reduction}).

The associated line bundle of the $\spinc$ structure
on $M \times \C$ is $\pr_1^* L \otimes \pr_2^* C$,
where $\pr_1$, $\pr_2$ are the projections from $M\times \C$
to $M$ and to $\C$, respectively.
Again by \S\ref{sec:spinc-reduction}, the line bundle associated
to the spin$^c$ structure $M_\cut$ is the reduced line bundle,
$$
	L_\cut = ( \pr_1^* L \otimes \pr_2^* C) _ \red \ ;
$$
reduced with respect to the anti-diagonal circle action.
On it we have a left circle action induced from the given action on $L$
and the trivial action on $\pr_2^* C$.

Over $M_+$ we have
\begin{equation} \labell{tensors}
 i_+^* L_\cut = \tiplus^* (\pr_1^* L \otimes \pr_2^* C)
 = L \otimes \C = L
\end{equation}
as $S^1$-equivariant complex line bundles,
because $\pr_1 \circ \tiplus = \operatorname{id}_{M_+}$,
and because the left circle action is trivial on the $\C$ piece.
This proves part 4(a).

\

Over $Z$ we have
$$\pi^* i_\red^* L_\cut
 = \tilde{i}_\red^* (\pr_1^* L \otimes \pr_2^* C)
 = i^* L \otimes \C$$
where $i$ is the inclusion of $Z$ into $M$,
and where $\C$ is the trivial complex line bundle
over $Z$ with the anti-standard circle action.  In particular,
$\C_\red = Z \times_{S^1} \C = N$, hence
$$ i_\red^* L_\cut = (i^* L \otimes \C)_\red = L_\red \otimes \C_\red
 = L_\red \otimes N.$$
This proves part 4(b),
thus completing the proof of Proposition \ref{cutting}.
\end{pf*}

In the course of the above proof we used the following lemma:

\begin{Lemma}
\labell{prodspinc}
The product of two equivariant spin$^c$ manifolds
has a natural equivariant $\spinc$ structure,
whose associated line bundle is the exterior tensor product of the
two original associated line bundles.
\end{Lemma}

\begin{pf}
The inclusion
$j: \SO(m) \times \SO(n) \hookrightarrow \SO(m+n)$
lifts to an inclusion
$j: \spin (m) \times_{\Z_2} \spin (n) \hookrightarrow \spin(m+n)$.
Define a group homomorphism
$ j : \spinc(m) \times \spinc(n) \to \spinc(m+n)$
by
$$ j( [A,a] , [B,b] ) =  [j(A,B),ab] $$
(where square brackets denote $\Z_2$-equivalence class in
$\spinc(k) = \spin(k) \times_{\Z_2} U(1)$).
The following diagrams commute:
\begin{equation}
\labell{fcommutes2}
   \begin{array}{ccc}
   \spinc(m) \times \spinc(n) & \stackrel{j}{\hookrightarrow} & \spinc(m+n) \\
   {\scriptstyle f \times f} \downarrow &  & {\scriptstyle f} \downarrow \\
   \SO(m) \times \SO(n) & \stackrel{j}{\hookrightarrow} & \SO(m + n)
   \end{array}
\end{equation}
and
\begin{equation}
\labell{gcommutes2}
   \begin{array}{ccc}
   \spinc(m) \times \spinc(n) & \stackrel{j}{\hookrightarrow} & \spinc(m + n)
\\
   {\scriptstyle \det \times \det} \downarrow
   \phantom{\scriptstyle \det \times \det}
	&  & {\scriptstyle \det} \downarrow \phantom{\scriptstyle \det}\\
   U(1) \times U(1)& \stackrel{\times}{\rightarrow} & U(1)
   \end{array}
\end{equation}
where the bottom arrow is multiplication,
and where $f$ and $\det$ are as in \S\ref{spinc quantization}.

Let $(M,P,p)$ and $(M',P',p')$ be spin$^c$ manifolds.
Let $\tilde{P}$ denote the
 $\spinc(n+m)$ bundle over $M \times M'$ that is associated
to the bundle $P \times P'$ by $j: \spinc(n) \times \spinc(m)
\to \spinc(n+m)$.
By diagram \eqref{fcommutes2} it is  clear that
$p \times p'$ naturally induces a $\spinc(n+m)$-equivariant
map $\tilde{p} : \tilde{P} \to \cF(T(M \times M'))$.
The final claim follows from diagram \eqref{gcommutes2}.
\end{pf}


\section{Proof of the main theorem}
\labell{proof-spinc}

We are now ready to prove our main theorem, Theorem~\ref{theorem-spinc},
that ``spin$^c$ quantization commutes with reduction''.
Following \cite{4}, to prove this theorem we apply the Atiyah-Segal-Singer
index theorem to the manifold $M$ and to the cut space $M_\cut$,
and then compare the results.

Let the circle group act on a compact manifold $M$.
Let $(P,p)$ be an equivariant spin$^c$ structure on $M$
with an associated line bundle $L$.
By the Atiyah-Segal-Singer index theorem \cite{ASIII},
the character, $\chi$, (defined in \S\ref{sec:intro})
of the equivariant index of the spin$^c$-Dirac
operator is given by a sum over the connected components
of the fixed point set:
\begin{equation} \labell{sum}
	\chi(\exp(u)) = \sum_F AS_F(u)\ , \quad u \in \Lie(S^1)\ .
\end{equation}
The contribution of the fixed component $F$ is:
\begin{equation}
\labell{expression}
	\AS_F (u) = (-1)^m \int_F \exp(\half \tc_1(L|_F)) \Aroof(F)
		   {{\Aroof}(NF) \over \Euler(NF)}\ ,
\end{equation}
where $L$ is the line bundle associated to the $\spinc$ structure,
where $NF$ is the normal bundle of $F$ in $M$ and has fiber dimension $2m$,
and where $\tc_1$, $\Aroof$, and $\Euler$ are, respectively,
the equivariant Chern class, A-roof class, and Euler class.
The orientations are chosen so that the orientation of $F$
followed by that of $NF$ gives the orientation of $M$.

Before proceeding with the proof, let us say a word about the derivation
of these formulas.
Formula \eqref{sum}, with Expression \eqref{expression},
follows from Theorem (3.9) in \cite{ASIII} by a derivation similar
to that of Formula (5.4) in \cite{ASIII}.
The difference is that Atiyah and Singer work out the case
of the Dirac operator associated to a $\spin$ structure,
whereas we need the Dirac operator associated to a $\spinc$ structure.
We now explain, briefly, the necessary adjustments:

The weight lattice of the maximal torus in the group
$\spinc(2n) = \spin(2n) \times_{\Z_2} U(1)$
is generated by the elements $x_1, \ldots, x_n, y$,
and $\half(x_1 + \ldots + x_n + y)$,
where $x_i$ and $y$ generate the weight lattices
of the maximal tori in $SO(2n) = \spin(2n) / \Z_2$
and in $U(1) / \Z_2$.
The $\spinc$-Dirac operator is constructed out of the
spinor representations
$\Delta^+_\C := \Delta^+ \otimes \C$ and
$\Delta^-_\C := \Delta^- \otimes \C$ of the group $\spinc$.
Temporarily adopting the conventions of Atiyah and Singer,
the weights occurring in these representations are
$\half (\pm x_1 \ldots \pm x_n + y)$
with an even number of minus signs for weights in $\Delta^+_\C$
and an odd number of minus signs for weights in $\Delta^-_\C$.
This gives $\ch(\Delta^+_\C) - \ch(\Delta^-_\C)
 = e^{y/2} \prod_{i=1}^n (e^{x_i/2} - e^{-x_i/2})$
instead of formula (5.1) in \cite{ASIII}.
The formula at the top of \cite[p.~571]{ASIII}
then gets multiplied by $e^{y/2}$,
giving $\text{index } D^+ = (-1)^n \exp(\half c_1(L)) \Aroof(M).$
Finally, since our operator agrees with theirs if $n = \half \dim M$
is even and disagrees if $n$ is odd
(see Remark \ref{global sign}), the factor $(-1)^n$ disappears.

We also note that in Formula (5.4) of \cite{ASIII} the sign depends ``in a
subtle global manner on the particular component" \cite[p.~572]{ASIII}.
This difficulty arises from the possible ambiguity in lifting
a $G$-action from the manifold to its spin bundle.
Since our group $G$ is connected, no such ambiguity arises.

We will now review the definitions of the characteristic classes
that occur in Expression \eqref{expression}.

First, if $L$ is a complex line bundle with ordinary Chern class
$[\omega]$, and if the circle group acts on the fibers of $L$ with
a weight $\alpha$, the equivariant Chern class of $L$
is $\omega + \alpha u$, where $u$ is a formal variable.
See Appendix~\ref{conventions}.

The Lie algebra of the maximal torus of $\SO(2n)$ is $\R^n$.
The Weyl group acts by permuting the $x_i$'s and by
replacing $x_i,x_j$ by $-x_i,-x_j$ for pairs of indices $i,j$.
Any power series on $\R^n$ that is invariant under these operations
determines a characteristic class on $\SO(2n)$-principal bundles.
In particular, the $\Aroof$-class is determined
from the power series expansion of the expression
\begin{equation}
\labell{Aroof}
	\Aroof = \prod_{i=1}^n {x_i \over e^{x_i/2}-e^{-x_i/2}}\ ,
\end{equation}
the Euler class is determined by the polynomial
\begin{equation}
\labell{Euler}
	\Euler = \prod_{i=1}^n x_i\ ,
\end{equation}
and their quotient is determined by the power series expansion of
the expression
\begin{equation}
\labell{A over e}
	\prod {1 \over e^{x_i/2} - e^{-x_i/2}}\ .
\end{equation}

In \eqref{expression}, $\Aroof(F)$ means the $\Aroof$-class associated
to the frame bundle $\cF(TF)$, on which the circle action is trivial.
For the normal bundle, $NF$, we need to work with the equivariant analogues
of the above characteristic classes.

Since the circle acts on the fibers of $NF$ while fixing
only the zero section, $NF$ is isomorphic
(as an equivariant oriented real vector bundle)
to a complex vector bundle.
By the (equivariant) splitting principle we can pretend that $NF$
splits into a direct sum of complex line bundles, $NF = \bigoplus_i L_i$.
The Chern-Weil recipe then dictates that in \eqref{A over e}
we replace each $x_i$ by the equivariant Chern class,
$\omega_i + \alpha_i u$, of $L_i$, where the $\omega_i$
and $\alpha_i$ are the ordinary Chern class and the circle weight on $L_i$.
This yields
\begin{equation}
\labell{product}
	{\Aroof( NF) \over \text{Euler} (NF)} (u) =
	\prod_i {1 \over e^{ (\omega_i + \alpha_i u)/2 }
	               - e^{-(\omega_i + \alpha_i u)/2 }}.
\end{equation}

The sum $\sum_F \AS_F(u)$ is a polynomial in $z = \exp(u)$ and
$z^{-1} = \exp(-u)$,
because it is the character of the virtual representation $Q(M,P,p)$.
The (virtual) dimension of the fixed space, $\dim Q(M,P,p)^{S^1}$,
is equal to the coefficient of $1$ in this polynomial.

The individual contribution, $\AS_F(u)$, while not  a polynomial
in $z$ and $z^{-1}$, does extend to a meromorphic function
of $z \in \C \cup \{ \infty \}$, and can be expanded into
either a formal power series in $z^{-1}$ plus a polynomial in $z$,
or a formal power series in $z$ plus a polynomial in $z^{-1}$.
In either case, the dimension of the space of fixed vectors in $Q(M,P,p)$
is given by summing  the coefficients of $1$ over all $F$'s.

We now expand \eqref{expression} into a formal power series in $z\inv$.
If $\alpha_i > 0$, the $i$'th factor in \eqref{product} expands as
$$
	e^{-(\omega_i + \alpha_i u)/2}
	(1 + cz^{-1} + c' z^{-2} + \ldots).
$$
(Here, and throughout the rest of this paper,
$c,c',\ldots$ will serve as place holders for
arbitrary cohomology classes; {\em a priori}, the $c$'s in this formula bear
no relation to any $c$'s appearing later.)
If $\alpha_i < 0$, the $i$'th factor in \eqref{product} expands as
$$
	- e^{(\omega_i + \alpha_i u)/2}
	(1 + cz^{-1} + c' z^{-2} + \ldots).
$$

Plugging into \eqref{expression}, we get:
\begin{equation} \labell{we-get}
	AS_F(u) = (-1)^m
	z^{{\half (\alpha - \sum |\alpha_i| )}}
	\int_F c   + c'z^{-1} +c''z^{-2} + \ldots \  .
\end{equation}
A similar computation, when we expand in powers of $z$, gives:
\begin{equation} \labell{we-get2}
	AS_F(u) = (-1)^m
	z^{\half(\alpha + \sum |\alpha_i| )}
	\int_F c + c'z + c''z^2 + \ldots
\end{equation}

The lemma below follows immediately.

\begin{Lemma}
\labell{contributes}
Let $S^1$ act on a smooth manifold $M$.
Let $(P,p)$ be a an equivariant spin$^c$ structure on $M$.
\begin{enumerate}
\item
When we expand the Atiyah-Segal-Singer index formula \eqref{sum},
\eqref{expression} for the quantization $Q(M,P,p)$ as a formal power series
in $z$, the fixed point component $F$ does not contribute to
$\dim Q(M,P,p)^{S^1}$ if
$$
   \mbox{ fiber weight at } F 	>
   - \sum | \mbox{tangent weights at } F|\ .
$$
\item
When we expand the Atiyah-Segal-Singer index formula
for $Q(M,P,p)$ as a formal power series in $z^{-1}$,
the fixed point component $F$ does not contribute to $\dim Q(M,P,p)^{S^1}$
if
$$
   \mbox{ fiber weight at } F <
   \sum | \mbox{tangent weights at } F|\ .
$$
\end{enumerate}
\end{Lemma}

We are now ready to prove our main theorem.
Let $S^1$ act on a manifold $M$, and let $(P,p)$ be an equivariant
spin$^c$ structure.  Let $Z$ be a reducible hypersurface that splits $M$ into
$M_+$ and $M_-$. Assume that conditions \eqref{conditions-spinc} are satisfied.
Equivalently, assume:
\begin{equation} \labell{conditions-spinc-again}
   \begin{array}{lcl}
   F \subset M_+ & \Rightarrow &
   \mbox{ fiber weight at } F 	>
   - \sum | \mbox{tangent weights at } F| \\
   F \subset M_- & \Rightarrow &
   \mbox{ fiber weight at } F 	<
   \phantom{-} \sum | \mbox{tangent weights at }F|
   \end{array}
\end{equation}

We begin by considering the quantization of $M$.
By Lemma~\ref{contributes}, fixed points in $M_-$ contribute nothing
to $Q(M,P,p)^{S^1}$ when we take the expansion in $z^{-1}$; that is,
\begin{equation}
\labell{only M+}
	\dim Q(M,P,p) ^{S^1} = \sum_{F \subset M_+}
	\begin{array}{ll}
	\mbox{ coefficient of 1 in the expansion of $AS_F(u)$} \\
	\mbox{ as a formal power series in $z\inv$.}
	\end{array}
\end{equation}

Now consider the cut space, $M_\cut$.
By Proposition~\ref{cutting}, part 1,
the set of fixed points in $M_\cut$ is the union
$$
	M_\cut^{S^1} = i_+ (M_+^{S^1}) \sqcup i_\red ( M_\red),
$$
where $M_+^{S^1}$ is the set of fixed points in $M_+$.
Let us denote by $F'$ and $X'$, respectively, the images
in $M_\cut$ of connected components $F$ of $M_+^{S^1}$
and $X$ of $M_\red$. Then
\begin{equation}
\labell{M cut}
	\begin{array}{l}
	\dim Q(M_\cut,P_\cut,p_\cut)^{S^1}  = \\
	\\
	\phantom{+}
	\displaystyle{\sum_{F \subset M_+}}
	\begin{array}{ll}
	\mbox{ coefficient of 1 in the expansion of $AS_{F'}(u)$ }\\
	\mbox{ as a formal power series in $z\inv$ }
	\end{array} \\
	 +  \displaystyle{\sum_{X \subset M_\red}}
	\begin{array}{ll}
	\mbox{ coefficient of 1 in the expansion of $AS_{X'}(u)$}\\
	\mbox{ as a formal power series in $z\inv$.}
	\end{array}
	\end{array}
\end{equation}

Let $F$ be a connected component of $M_+^{S^1}$, and let
$F' := i_+(F)$ be its image in $M_\cut$.
By Proposition \ref{cutting}, parts 1, 3(a), and 4(a),
and by \eqref{expression},
the contribution of $F$ to the quantization of $M$
and the contribution of $F' = i_+(F)$ to the quantization of $M_\cut$
are identical:
\begin{equation}
\labell{equal1}
	AS_{F'} (u) = AS_F (u) \ .
\end{equation}

Let $X$ be a connected component of $M_\red$, let $X' := i_\red(X)$
be its image in $M_\cut$, and let $NX'$ be the normal bundle of $X'$
in $M_\cut$.
Let $N := Z \times_{S^1} \C$ be the line bundle over $X$ as in \eqref{normal}.
The contribution \eqref{expression} of the fixed component $X'$
to the quantization of $M_\cut$ is:
\begin{eqnarray}
 AS_{X'} (u) & = &
	- \int_{X'} \exp(\half \tc_1(L_\cut|_{X'})) \Aroof(X')
	{\Aroof(NX') \over \text{Euler}(NX')} \notag \\
 & = &  - \int_X \exp (\half \tc_1 (L_\red \otimes N) )
	\Aroof(X) {\Aroof(N) \over \text{Euler}(N)} \notag \\
 & = &  - \int_X \exp (\half (c_1 (L_\red) ) )
	\Aroof(X) {\exp (\half \tc_1 (N)) \Aroof(N)
	\over \text{Euler}(N)} \ ,		\labell{work}
\end{eqnarray}
where the middle equality follows from Proposition \ref{cutting},
parts 2, 3(b), and 4(b).

Since $N$ is a complex line bundle with fiber weight $1$,
and by \eqref{product}, we have
\begin{eqnarray}
\exp(\half \tc_1(N)) & = & \exp(\half(\omega+u))
\quad , \quad \mbox{and}  \\
{\Aroof(N) \over \Euler(N)}  & = &
{1 \over \exp(\half(\omega+u)) - \exp(-\half(\omega+u)) } \ ,
\end{eqnarray}
where $[\omega]$ is the first Chern class of $N$.
Plugging this into \eqref{work} and expanding in powers of
$z\inv = \exp(-u)$, we get
$$
 AS_{X'} (u) = - \int_X \exp (\half (c_1 (L_\red) ) )
	\Aroof(X) (1 + cz^{-1} + c' z^{-2} + \ldots )\ .
$$
By applying the (non-equivariant) Atiyah-Singer formula
to $M_\red$, we see that
\begin{equation}
\labell{equal2}
\begin{array}{l}
	\displaystyle{\sum_{X \subset M_\red}}
	\begin{array}{ll}
	\mbox{ coefficient of 1 in the expansion of $AS_{X'}(u)$} \\
	\mbox{ as a formal power series in $z\inv$}
	\end{array}  \hfill \\
\hfill  = -\dim Q(M_\red,P_\red,p_\red)\ .
\end{array}
\end{equation}

Comparing the expressions \eqref{only M+}, \eqref{M cut}, \eqref{equal1}
and \eqref{equal2}, we get:
\begin{equation}
\labell{left}
	\dim Q(M_\cut,P_\cut,p_\cut)^{S^1} =
	\dim Q(M,P,p)^{S^1} - \dim Q(M_\red,P_\red,p_\red)\ .
\end{equation}

Now consider the expansions in $z$.
For every fixed component of the form $F' = i_+(F)$,
for $F \subset M_+$,
the first condition of \eqref{conditions-spinc-again}
together with part 1 of Lemma~\ref{contributes}
implies that this fixed point set does not contribute.
Similarly, consider a fixed component of the form $X' = i_\red(X)$,
for $X \subset M_\red$.
Its fiber weight is $\alpha = 1$,
by Proposition \ref{cutting}, part 2.
Hence, by the part 1 of Lemma~\ref{contributes}, $X'$ does not
contribute as well. Therefore,
\begin{equation}
\labell{right}
	\dim  Q(M_\cut, P_\cut, p_\cut) ^{S^1} = 0\ .
\end{equation}

Since \eqref{left} and \eqref{right} must be equal, we have
$$
	\dim Q(M,P,p)^{S^1} = \dim Q(M_\red,P_\red,p_\red)\ ,
$$
{\em i.e.\/}, spin$^c$ quantization commutes with reduction.


\appendix

\section{Remark on conventions}
\labell{conventions}

Let $\pi: P \to M$ be a principal $S^1$ bundle.
A connection form on $P$ is an invariant one-form,
$\Theta \in \Omega^1(P,\Lie(S^1))$,
that restricts to the canonical form on each fiber.
Its curvature is the two-form
$\Omega \in \Omega^2(M,\Lie(S^1))$ such  that
$\pi^*(\Omega) = d\Theta$.

If $S^1$ acts on $P$,
define a map $\Psi : M  \to (\Lie(S^1))^* \otimes \Lie(S^1)$
by $\pi^* \Psi (\xi)= \Theta(\xi_P)$,
where $\xi_P$ is the vector field on $P$
corresponding to $\xi \in \Lie(S^1)$ under the action.
Then $\Psi$ is a moment map for $\Omega$, that is,
$d \Psi (\xi) = -\iota(\xi_M) \Omega$ for all $\xi \in \Lie(S^1)$,
where $\xi_M$ is the vector field on $M$ that generates the action
of $\exp(t\xi)$, $t \in \R$.

Over each component $F$ of the fixed point set,
this action gives a representation $\rho_F : S^1 \to S^1$.
The value of $\Psi$ at $F$ is precisely the weight
of the representation $\rho_F$; no choices are necessary.
Under the natural identification of $\Lie(S^1)^* \otimes \Lie(S^1)$
with $\R$, the weights are identified with the integers, $\Z \subset \R$,
and $\Psi$ is identified  with a real valued function $\Phi : M \to \R$.

In this paper, we always identify $\Lie(S^1)$ with $\R$
by the isomorphism that takes the positive primitive
lattice element, $\eta \in \Lie(S^1)$, to $-1$.
The curvature, $\Omega$, is then identified with a real-valued two-form,
$\omega$, which represents the first Chern class of $P$.
The criterion for a moment map then becomes $d \Phi = \iota (\eta_M) \omega$.


\end{document}